\documentclass[a4paper]{article}
\usepackage[francais]{babel}
\usepackage[latin1]{inputenc}
\usepackage{amsmath}
\usepackage{array}
\usepackage{amsthm}
\usepackage{amssymb}
\input xy
\xyoption{all}

\author{Aur\'elien DJAMENT}

\title{Le foncteur $V\mapsto\FF[V]^{\otimes 3}$ entre $\FF$-espaces
  vectoriels est noethérien}

\date{Septembre 2008 (première version : février 2007)}
\newcommand{\A}{{\mathcal{A}}}

\newcommand{\C}{{\mathcal{C}}}
\newcommand{\D}{{\mathcal{D}}}
\newcommand{\F}{{\mathcal{F}}}
\newcommand{\E}{{\mathcal{E}}}
\newcommand{\R}{{\mathcal{R}}}
\newcommand{\Prl}{{\mathcal{P}}}

\newcommand{\K}{{\mathcal{K}}}
\newcommand{\Q}{{\mathcal{Q}}}
\newcommand{\N}{{\mathcal{N}}}
\newcommand{\FF}{{\mathbb{F}_2}}
\newcommand{\p}{{\mathfrak{p}}}

\newcommand{\Gr}{{\mathcal{G}r}}

\newcommand{\T}{{\mathcal{T}}}

\newtheorem{thm-intro}{Théorème}
\newtheorem{conj-intro}[thm-intro]{Conjecture}

\newtheorem{theo}{Théorème}[section]
\newtheorem{pr}[theo]{Proposition}
\newtheorem{cor}[theo]{Corollaire}
\newtheorem{lm}[theo]{Lemme}

\newtheorem{conj}[theo]{Conjecture}

\theoremstyle{definition}
\newtheorem{defi}[theo]{Définition}
\newtheorem{nota}[theo]{Notation}
\newtheorem{hyp}[theo]{Hypothèse}

\newtheorem{def-intro}[thm-intro]{Définition}

\theoremstyle{remark}
\newtheorem{rem}[theo]{Remarque}
\newtheorem{ex}[theo]{Exemple}

\begin{document}

\maketitle

\begin{abstract} Les foncteurs entre espaces vectoriels, ou {\em représentations génériques} des groupes linéaires (\cite{K1}, \cite{K2}), interviennent en topologie algébrique et en $K$-théorie comme en théorie des représentations (\cite{FFPS}). Nous présentons ici une nouvelle méthode pour aborder les problèmes de finitude et la dimension de Krull dans ce contexte.

Plus précisément, nous démontrons que, dans la catégorie $\F$ des foncteurs
entre espaces vectoriels sur $\FF$, le produit tensoriel entre
$P^{\otimes 3}$, où $P$ désigne le
 foncteur projectif $V\mapsto\mathbb{F}_2[V]$, et un foncteur de longueur finie est noethérien et déterminons sa structure. Seul était antérieurement connu le caractère
noethérien de $P^{\otimes 2}\otimes F$ pour $F$ de longueur finie.

Nous utilisons pour cela plusieurs foncteurs de division, dont nous
analysons l'effet sur des foncteurs de type fini de $\F$ à l'aide des
catégories de foncteurs en grassmanniennes, introduites dans \cite{artsmf}. Cela nous permet de ramener le problème initial à des calculs explicites finis portant sur des représentations modulaires de
groupes linéaires (où intervient notamment la représentation de Steinberg), qui renseignent finalement sur des phénomènes infinis en théorie des représentations. 

\begin{center}{\bf Abstract}
\end{center}

Functors between vector spaces, or {\em generic representations} of linear groups (\cite{K1}, \cite{K2}  intervene in algebraic topology and in $K$-theory as in representation theory (cf. \cite{FFPS}). We present here a new method to approach finiteness problems and the Krull dimension in this context.

More precisely, we prove that, in the category $\F$ of functors
  between $\FF$-vector spaces, the tensor product between $P^{\otimes
    3}$, where $P$ denotes the projective functor $V\mapsto\mathbb{F}_2[V]$,
    and a functor of finite length is noetherian and we determine its structure. The only case known to date was the noetherian character
  of $P^{\otimes 2}\otimes F$ for $F$ of finite length. 
 
For this we use several division functors, whose effect on finitely
generated functors  of $\F$ is analyzed with the help of
grassmannian functor categories, introduced in \cite{artsmf}. It allows us to reduce the initial problem to finite calculations on modular representations of linear groups (where appears in particular Steinberg's representation), which inform ultimately on infinite phenomena in representation theory.
\end{abstract}

\smallskip

\noindent
\begin{small}{\em Mots-clés : } cat\'egories de foncteurs, repr\'esentations
modulaires, groupes linéaires, foncteurs de division, filtration de
Krull, grassmanniennes.

\noindent
{\em Classification math. : } 18A25 (16P60 18D10 18E25 20C33 55S10).
\end{small}

\tableofcontents

\section*{Introduction}

La catégorie des {\em représentations génériques} des groupes
linéaires sur le corps à deux éléments $\FF$ désigne, depuis les
articles de
Kuhn \cite{K1} et \cite{K2}, la
catégorie $\mathbf{Fct}(\E^f,\E)$, notée $\F$, des foncteurs de $\E^f$ vers $\E$, où l'on note
$\E$ la catégorie des $\FF$-espaces vectoriels et $\E^f$ celle des
espaces vectoriels de dimension finie --- le corps de base $\FF$ sera par
la suite implicite. \'Etudiée par les topologues, en raison de
ses liens avec l'algèbre de Steenrod, cette catégorie abélienne possède aussi des
liens cohomologiques profonds avec les groupes linéaires ; nous renvoyons pour cela le
lecteur à l'ouvrage de synthèse \cite{FFPS}. La catégorie $\F$ possède
des objets simples classifiés à l'aide de la théorie des
représentations des groupes symétriques ou linéaires, et l'on dispose
d'outils de calcul cohomologique efficaces sur les foncteurs de
longueur finie. Cependant, la compréhension intrinsèque de la
catégorie $\F$ comme de ses liens avec les différents groupes
linéaires se complique énormément dès que l'on passe des foncteurs
simples (ou de
longueur finie) à des foncteurs de type fini. En effet, la structure de
nombreux foncteurs demeure inconnue, et la plupart des assertions
reliant les propriétés des facteurs de composition d'un foncteur à sa
structure globale restent conjecturales.

\smallskip

Pour illustrer ces difficultés, rappelons que la catégorie $\F$ est engendrée par les foncteurs projectifs de type
fini $P_V$, dits {\em standard}, définis par $P_V=\FF[{\rm
  hom}_{\E^f}(V,-)]$. Ici $V$ désigne un objet de $\E^f$ et $\FF[E]$ l'espace vectoriel librement
engendré par un ensemble $E$ ; le lemme de Yoneda fournit un
isomorphisme naturel ${\rm hom}_\F(P_V,F)\simeq F(V)$. L'existence de
ces générateurs montre que $\F$ est une {\em catégorie de Grothendieck} (i.e. une catégorie abélienne possédant des colimites filtrantes exactes et un ensemble de générateurs\,\footnote{Le lecteur est invité à se reporter à l'ouvrage de
  Popescu \cite{Pop} pour les généralités sur les catégories
  abéliennes.}).

On notera $P$ le foncteur $P_\FF$, lorsqu'aucune confusion ne peut en résulter ;
remarquer que l'on a $P_V\simeq P^{\otimes \dim V}$, le produit
tensoriel de $\F$ étant calculé au but. Les problèmes soulevés par
l'étude des foncteurs projectifs standard sont illustrés par la
conjecture suivante, ouverte depuis une quinzaine d'années.

\begin{conj-intro}\label{ca} Pour tout $n\in\mathbb{N}$, le foncteur
  $P^{\otimes n}$ est noethérien (i.e. toute suite croissante de
  sous-foncteurs de $P^{\otimes n}$ stationne).
\end{conj-intro}

Cet énoncé, équivalent au caractère localement noethérien de la
catégorie $\F$, est connu sous le nom de {\em conjecture artinienne} car sa
formulation initiale, émise par Lannes et Schwartz à partir du
contexte topologique sous-jacent à $\F$, était donnée en termes duaux
d'objets artiniens. Rappelons ici que la catégorie $\F$ possède un
foncteur de dualité $D : \F^{op}\to\F$ tel que $D(F)(V)=F(V^*)^*$, où
$(-)^*$ désigne le foncteur de dualité usuel des espaces vectoriels ;
la plupart des premières investigations sur la structure de $\F$ ont
été menées sur les objets injectifs $I_V=DP_V$.

La conjecture artinienne~\ref{ca}, discutée en détails dans les
articles \cite{GP1}
et \cite{artsmf}, est facile pour $n\leq 1$ ($P_0=\FF$, foncteur constant, est simple ; $P$ est somme directe de $\FF$ et d'un foncteur unisériel de facteurs de composition les puissances extérieures) ; elle a été démontrée par
Powell pour $n=2$ dans \cite{GP5}. Nous avons étendu le caractère
noethérien de $P^{\otimes 2}$ à $P^{\otimes 2}\otimes F$, où $F\in {\rm
  Ob}\,\F$ est un foncteur de longueur finie, dans l'article
\cite{artsmf}, grâce à l'introduction des {\em catégories de
  foncteurs en grassmanniennes}. Ces
catégories de foncteurs, dont la définition et quelques propriétés sont rappelées dans la
section~\ref{s1}, permettent non seulement de progresser dans l'étude
de la conjecture artinienne, mais aussi de préciser nettement son
énoncé (forme~\guillemotleft~extrêmement forte~\guillemotright), à l'aide d'une description conjecturale de la filtration de
Krull de la catégorie $\F$ (cf. conjecture~\ref{caef2}). Le travail du
dernier chapitre de
\cite{artsmf} ramène la conjecture artinienne dans sa
forme renforcée à celui de la description des foncteurs nilpotents
pour des endofoncteurs de $\F$, notés $\tilde{\nabla}_n$, introduits par Powell dans \cite{GP4},
qu'il a lourdement utilisés pour démontrer son résultat relatif à
$P^{\otimes 2}$ --- le théorème principal dudit chapitre de \cite{artsmf} étend
d'ailleurs, dans le contexte des catégories de foncteurs en
grassmanniennes, le théorème de structure donné par Powell dans
\cite{GP2}. 

\medskip

Dans \cite{GP5}, les renseignements sur le foncteur $\tilde{\nabla}_2$
nécessaires à établir la conjecture artinienne pour $n=2$ sont obtenus
par une méthode assez directe, inspirée par la théorie des
représentations. De fait, Powell remarque que les facteurs de
composition possibles d'un foncteur $F$ tel que
$\tilde{\nabla}_2(F)=0$ sont limités --- ce sont forcément des puissances
extérieures ---  et que, de plus, certaines extensions entre
puissances extérieures ne peuvent apparaître comme sous-quotient d'un tel
foncteur. Il s'appuie alors sur le calcul des extensions entre
puissances extérieures, qui s'avèrent peu nombreuses, pour
conclure. La généralisation de cette approche se heurte au problème
majeur suivant : pour $n\geq 3$, les extensions entre les foncteurs
simples annihilés par $\tilde{\nabla}_n$ sont beaucoup plus nombreuses
qu'entre des puissances extérieures et elles ne sont de plus pas
comprises en général.

Dans le cas $n=3$, le présent article contourne cet obstacle ; nous
parvenons à estimer les foncteurs annihilés par $\tilde{\nabla}_3$,
aboutissant au théorème suivant.

\begin{thm-intro}\label{thi} Si $F$ est un foncteur de longueur finie de $\F$, le
  foncteur $P^{\otimes 3}\otimes F$ est noethérien. En particulier, $P^{\otimes 3}$ est noethérien.
\end{thm-intro}

Notre observation initiale consiste à remarquer que le foncteur
$\tilde{\nabla}_3$ peut être avantageusement remplacé, pour les
considérations d'annulation qui nous intéressent, par certains
{\em foncteurs de division}, c'est-à-dire par des adjoints à gauche au
produit tensoriel par un foncteur. D'une manière générale, pour tout
entier $n\geq 1$, deux foncteurs de division, dont l'un est exact,
annihilent exactement les mêmes foncteurs simples de $\F$ que
$\tilde{\nabla}_n$. Les deux foncteurs
en question proviennent de la {\em représentation de Steinberg} de
$GL_n(\FF)$, dont nous avons déjà observé, dans la section~16.1 de \cite{artsmf}, qu'elle
peut jouer un rôle important dans l'étude de la structure de la catégorie $\F$.

Notre manipulation des foncteurs de
division~\guillemotleft~remplaçant~\guillemotright~le foncteur
$\tilde{\nabla}_3$ repose sur plusieurs idées. D'une part, nous
n'abordons pas directement le problème de l'annulation de ces foncteurs
de division, mais celui de foncteurs de division plus simples, que
nous utilisons ensuite abondamment : heuristiquement, nous cherchons à
produire successivement des foncteurs $A$ tels que l'on contrôle la taille d'un foncteur
raisonnable $F$ dès que la taille de la division de $F$ par $A$ est contrôlée.

Pour mener à bien ce programme, nous établissons des propriétés
générales des foncteurs de division qui reposent sur les catégories de
foncteurs en grassmanniennes. On peut résumer, de manière imprécise,
les résultats de la manière suivante : 
\begin{enumerate}\item les foncteurs de division des catégories de
  foncteurs en grassmanniennes préservent les foncteurs de longueur
  finie et diminuent (au sens large) leur taille ; 
\item on contrôle le terme principal de la division d'un foncteur en grassmanniennes de
  longueur finie $X$ ; si $X$ est assez gros, ce terme principal est
  du même ordre de grandeur que $X$ ;
\item en revanche, les foncteurs dérivés de degré strictement positif
  des foncteurs de division diminuent toujours strictement la taille
  des foncteurs en grassmanniennes de longueur finie.
\end{enumerate}
La dernière propriété s'avère essentielle pour pallier l'inexactitude
des foncteurs de division et permettre certains raisonnements inductifs utilisant les suites exactes longues d'homologie associées aux foncteurs de division pour aborder la conjecture~\ref{ca}.

Outre ces propriétés essentiellement formelles, nos arguments
nécessitent des ingrédients propres à assurer le fonctionnement
d'arguments de récurrence emboîtant différents foncteurs de
division. Ceux-ci proviennent de la théorie des représentations : des calculs concrets de produits tensoriels de représentations
simples des groupes linéaires interviennent, dans des cas particulièrement
élémentaires du point de vue des représentations, ce qui explique que la
généralisation de notre approche au cas $n>3$ de la conjecture
artinienne~\ref{ca} pose des problèmes encore non élucidés.

\medskip

Dans l'article \cite{art1}, nous utilisions déjà des foncteurs de
division dans la catégorie $\F$ pour progresser dans l'étude de la
conjecture artinienne. Il s'agissait, comme dans le présent travail,
de combiner des propriétés formelles de ces foncteurs à des
considérations concrètes issues de la théorie des
représentations. Cependant, les deux démarches présentent de notables
différences conceptuelles.

Tout d'abord, \cite{art1} adopte des conventions duales de
celles de cet article : celui-là s'intéresse aux injectifs $I_V$, alors que
nous travaillons dans celui-ci sur les projectifs $P_V$. Pour traduire en termes
d'objets projectifs la démarche de \cite{art1}, il faut donc
s'attacher aux foncteurs duaux des foncteurs de division, c'est-à-dire
les foncteurs hom internes (adjoints à {\em droite} au produit
tensoriel). De fait, les comportements des foncteurs hom internes et
des foncteurs de division diffèrent profondément sur les foncteurs
$P_V$ : si $F$ est un foncteur non nul, la division de $P_V$ par $F$
est toujours non nulle pourvu que la dimension de $V$ soit assez
grande, tandis que l'image de $P_V$ par le foncteur hom interne de
source $F$ est toujours nulle si $F$ est de longueur finie et sans
terme constant. Pour autant, les foncteurs hom interne
et de division se comportent de manière analogue sur les
foncteurs simples de $\F$, qui sont auto-duaux, de sorte que certains
raisonnements sur les facteurs de composition employés dans cet article et
dans \cite{art1} se rejoignent partiellement.

Une autre différence fondamentale entre les deux approches que nous
discutons réside dans le principe même d'utilisation des foncteurs de
division pour comprendre la structure d'un objet de la catégorie
$\F$. Dans \cite{art1}, toute la stratégie repose sur la détection
de facteurs de composition significatifs dans certains foncteurs :
il s'agit d'une stratégie~\guillemotleft~locale~\guillemotright, qui étudie
les foncteurs à partir de leurs constituants élémentaires. Ici, notre
stratégie est~\guillemotleft~globale~\guillemotright~en ce sens
qu'elle consiste au contraire à étudier les objets de
$\F$ en leur appliquant des foncteurs de
division, ce qui en général a pour effet, entre autres, d'en annihiler de nombreux
facteurs de composition. Elle pourrait d'ailleurs souvent s'exprimer
en termes de catégories quotients, à ceci près qu'on ne peut procéder
ainsi sans de lourdes précautions en raison de l'inexactitude des
foncteurs de division (c'est là qu'intervient le contrôle sur les
foncteurs dérivés évoqué précédemment).

\smallskip

Pour terminer les remarques générales relatives aux techniques
développées dans le présent travail, soulignons les nombreux détours
nécessaires pour parvenir au résultat principal. La propriété des
foncteurs annihilés par $\tilde{\nabla}_3$ dont elle se déduit
requiert l'emploi successif de nombreux foncteurs auxiliaires et
contient plusieurs récurrences imbriquées. La manière hautement non
explicite dont elle se démontre illustre l'extrême complexité
combinatoire des sous-foncteurs de $P^{\otimes 3}$ (et a fortiori des
$P^{\otimes 3}\otimes F$). La situation semble encore plus intrigante
au-delà, car le cadre formel de cet article ne suffit pas à prouver le
caractère noethérien de $P^{\otimes 4}$. Les améliorations
envisageables pour y parvenir posent de grandes difficultés
relatives aux représentations des groupes linéaires qui font
apparaître le cas de $P^{\otimes 3}$ --- pourtant nettement plus ardu que
celui de $P^{\otimes 2}$, où des techniques encore assez explicites
fonctionnent --- comme miraculeusement élémentaire.

\medskip

Cet article s'organise comme suit. La première section introduit nos
notations générales et présente les rappels nécessaires sur les
catégories de foncteurs en grassmanniennes (ce
matériel se trouve dans \cite{artsmf}). La deuxième section
fournit tous les outils formels pour manipuler les foncteurs de
division dans l'optique de la conjecture artinienne (on pourra en première lecture se concentrer sur les énoncés
des propositions~\ref{pr-cadiv} et~\ref{pr-omd} puis aborder directement la dernière section). Enfin, la
troisième section donne la démonstration du théorème~\ref{thi}, à
partir d'énoncés généraux mettant en évidence les éléments non
formels issus de la théorie des représentations nécessaires pour la
mener à bien. L'appendice rappelle quant à lui quelques faits classiques sur la catégorie $\F$ et les notions de finitude afférentes ; il fixe les notations de paramétrisation des foncteurs simples et des représentations simples des groupes linéaires à l'aide des partitions.

\bigskip

L'auteur adresse ses chaleureux remerciements à Geoffrey Powell pour
ses discussions sur la catégorie~$\F$ et ses conseils qui ont
permis d'améliorer ce texte. Il témoigne aussi
sa gratitude à Lionel Schwartz pour ses conversations sur la conjecture artinienne et à Vincent Franjou pour ses utiles remarques sur la première version de ce travail.

\section{Les foncteurs en grassmanniennes}\label{s1}

Cette section expose sans démonstration les constructions et résultats principaux de l'article \cite{artsmf}, qui mène l'étude de la catégorie $\F$ à partir de catégories de foncteurs auxiliaires. Le lecteur peut se reporter à l'appendice pour ce qui concerne les définitions, notations et premiers faits classiques sur la catégorie~$\F$.

\paragraph*{Les catégories $\F_\Gr$ et $\F_{\Gr,n}$.} L'outil
fondamental introduit dans \cite{artsmf} pour étudier la structure de la
catégorie $\F$ est donné par les {\em catégories de foncteurs en
grassmanniennes}, dont les plus importantes sont les suivantes. On
désigne par $\E^f_\Gr$ la catégorie des couples $(V,W)$, où $V$ est un
objet de $\E^f$ et $W$ un sous-espace vectoriel de $V$ (on notera
$\Gr(V)$ l'ensemble de ces sous-espaces) ; les morphismes $(V,W)\to
(V',W')$ de $\E^f_\Gr$ sont les flèches $f : V\to V'$ de $\E^f$ telles
que $f(W)=W'$. Si $n$ est un entier, on introduit les sous-catégories pleines $\E^f_{\Gr,n}$ et $\E^f_{\Gr, \leq n}$ des objets $(V,W)$ tels que $\dim W=n$ et
$\dim W\leq n$ respectivement. On pose $\F_{\Gr}=\mathbf{Fct}(\E^f_{\Gr},\E)$, $\F_{\Gr,n}=\mathbf{Fct}(\E^f_{\Gr,n},\E)$ et $\F_{\Gr,\leq n}=\mathbf{Fct}(\E^f_{\Gr,\leq n},\E)$.
Les inclusions $\E^f_{\Gr,n}\hookrightarrow\E^f_{\Gr}$ et $\E^f_{\Gr,\leq n}\hookrightarrow\E^f_{\Gr}$ induisent par précomposition des foncteurs de restrictions
$\R_{n} : \F_{\Gr}\to\F_{\Gr,n}$ et $\R_{\leq n} : \F_{\Gr}\to\F_{\Gr,\leq n}$. On dispose également de foncteurs de {\em prolongement par zéro}
$\Prl_{n} : \F_{\Gr,n}\to\F_{\Gr}$ et $\Prl_{\leq n} : \F_{\Gr,\leq n}\to\F_{\Gr}$, à l'aide desquels on peut identifier $\F_{\Gr,n}$ et $\F_{\Gr,\leq n}$ à des sous-catégories
épaisses de $\F_{\Gr}$.

La précomposition par le foncteur d'oubli $\E^f_{\Gr}\to\E^f\quad (V,W)\mapsto V$ fournit un foncteur $\iota : \F\to\F_{\Gr}$. Il possède un adjoint à gauche $\omega : \F_{\Gr}\to\F$, appelé
foncteur d'{\em intégrale en grassmanniennes}. Explicitement, on a 
$$\omega(X)(V)=\bigoplus_{W\in\Gr(V)}X(V,W).$$
On note $\omega_{n}$ le foncteur composé $\F_{\Gr,n}\xrightarrow{\Prl_{n}}\F_{\Gr}\xrightarrow{\omega}\F$. Les foncteurs $\omega$ et $\omega_{n}$ sont exacts (ils commutent même à toutes les limites
et colimites).

\begin{rem} Les foncteurs $\omega$ et $\omega_n$ constituent l'outil essentiel de la théorie ; leur principale vertu est de transformer de petits foncteurs (par exemple, les foncteurs constants) en des foncteurs beaucoup plus gros. Il s'agit ainsi de ramener la compréhension d'objets infinis de $\F$ à celle d'objets finis des catégories de foncteurs en grassmanniennes.
\end{rem}

\paragraph*{Foncteurs différences et foncteurs finis.} Pour $E\in {\rm
  Ob}\,\E^f$, on introduit des endofoncteurs de $\F_\Gr$,
$\F_{\Gr,n}$ et $\F_{\Gr,\leq n}$ appelés {\em décalages par $E$},
tous notés  $\Delta_E^\Gr$, donnés par 
$\Delta^\Gr_E(X)(V,W)=X(V\oplus E,W)$ ; ils sont analogues au foncteur décalage $\Delta_E$ de la catégorie $\F$ (cf. appendice, page~\pageref{fdif}).

On introduit ensuite des endofoncteurs de $\F_\Gr$,
$\F_{\Gr,n}$ et $\F_{\Gr,\leq n}$ appelés {\em
  foncteurs différences} et notés $\Delta^\Gr$ par le scindement $\Delta_\FF^\Gr\simeq id\oplus\Delta^\Gr_\FF$. Les foncteurs décalages et différences
commutent à toutes les limites et colimites, comme dans le cas de la catégorie~$\F$. 

Les notions de {\em foncteurs polynomiaux} et de {\em degré} d'un foncteur polynomial se définissent dans $\F_\Gr$,
$\F_{\Gr,n}$ et $\F_{\Gr,\leq n}$ comme dans $\F$, en remplaçant $\Delta$ par $\Delta^\Gr$.

Nous utiliserons
l'analogue suivant de la proposition~\ref{rap-f} :
\begin{pr}[\cite{artsmf}, section~5.5]\label{pol-fgr}\begin{itemize}\item Dans les
    catégories $\F_{\Gr,n}$ et $\F_{\Gr,\leq n}$, les foncteurs finis
    sont les foncteurs polynomiaux à valeurs de dimension finie.
\item Un foncteur $X$ de $\F_\Gr$ est fini si et seulement s'il est
  polynomial, à valeurs de dimension finie et qu'il existe un entier
  $n$ tel que $X$ appartienne à l'image essentielle de $\Prl_{\leq n} : \F_{\Gr,\leq n}\to\F_{\Gr}$.
\end{itemize}
\end{pr}

Pour alléger, nous noterons simplement $\F_\Gr^f$ la catégorie
$(\F_\Gr)^f$ des objets finis de $\F_\Gr$ (cf. appendice) et adopterons d'autres
simplifications d'écriture analogues.

\begin{rem}\label{rqf} Dans la catégorie $\F$, les foncteurs polynomiaux de degré (au plus)
$0$ sont les foncteurs constants. Dans les catégories de foncteurs en
grassmanniennes, il y a plus de foncteurs de degré nul ; ces foncteurs
sont nommés {\em foncteurs pseudo-constants}. 
\end{rem}

Nous donnons maintenant leur
description dans le cas de $\F_{\Gr,n}$. 

\begin{nota}\label{dfpc}\begin{itemize}\item On note $\varepsilon_n :
\F_{\Gr,n}\to GL_n - \mathbf{mod}$ le foncteur d'évaluation donné sur les objets par $\varepsilon_n(X)=X(\FF^n,\FF^n)$.
\item Si $W$ est un espace vectoriel de dimension $n$, on désigne par ${\rm Iso}(\FF^n,W)$ le
$GL_n$-ensemble à droite libre et transitif des isomorphismes de
$\FF^n$ sur $W$.
\item On note $\rho_n : GL_n - \mathbf{mod}\to\F_{\Gr,n}$ le foncteur donné par $\rho_n(M)(V,W)=\FF[{\rm Iso}(\FF^n,W)]\underset{GL_n}{\otimes}
M$.
\end{itemize}
\end{nota}

\begin{pr}\label{prpc} \begin{itemize}\item Le foncteur $\rho_n$ est adjoint à gauche au foncteur $\varepsilon_n$.
\item Le noyau du foncteur différence de $\F_{\Gr,n}$ (i.e. la sous-catégorie pleine des foncteurs pseudo-constants) égale l'image essentielle du foncteur $\rho_n$.
\end{itemize}
\end{pr}

Les foncteurs définis ci-après sont les~\guillemotleft~constituants élémentaires~\guillemotright~des foncteurs pseudo-constants.

\begin{defi}\label{df-pow} Pour $\lambda\in\p_n$, on pose
$Q_\lambda=\omega_n\rho_n(R_\lambda)$. Un tel foncteur est appelé {\em
  foncteur de Powell}.
\end{defi}

Rappelons quelques faits utiles sur ces foncteurs (cf. \cite{GP2}, où
les duaux des $Q_\lambda$ sont étudiés sous le nom de {\em foncteurs
  co-Weyl}, et \cite{artsmf}, section~2.3, pour une présentation analogue à celle-ci) : 

\begin{pr}\label{prpw}\begin{itemize}\item Le cosocle\,\footnote{i.e. plus grand quotient semi-simple, cf. appendice.} de $Q_\lambda$ est $S_\lambda$ ;
\item si $X\in
{\rm Ob}\,\F^f_{\Gr,\leq n}$ est pseudo-constant, alors
$\omega\Prl_{\leq n}(X)$
possède une filtration finie dont les sous-quotients sont des
foncteurs de Powell $Q_\lambda$ avec $\lambda_1\leq n$. Cela vaut en
particulier pour le foncteur projectif $P^{\otimes n}$.
\end{itemize}
\end{pr}

\begin{nota}\label{nvnot} Soient $n$ et $k$ deux entiers.
On désigne par $\F_{\Gr}[n,k]$ la sous-catégorie pleine (qui est épaisse) de $\F_{\Gr}^f$ formée des foncteurs $X$ tels que :
\begin{enumerate}\item $X$ appartient à l'image essentielle du foncteur $\Prl_{\leq n} : \F_{\Gr,\leq n}\to\F_{\Gr}$ ;
\item le foncteur $\R_{n}(X)$ de $\F_{\Gr,n}$ est de degré inférieur à $k$.
\end{enumerate}\end{nota}
Si l'on munit $(\mathbb{N}\cup\{-1\})^2$ de l'ordre lexicographique,
on a $\F_{\Gr}[n,k]\subset\F_{\Gr}[n',k']$ pour $(n,k)\leq (n',k')$ ;
la proposition~\ref{pol-fgr} montre que $\F_\Gr^f$ est la réunion des
sous-catégories $\F_{\Gr}[n,k]$. Nombre des arguments de
récurrence utilisés dans la section~\ref{s-div} seront relatifs à ce bon ordre sur
$(\mathbb{N}\cup\{-1\})^2$, utilisé dans ce
contexte. C'est ce qui motive l'introduction de la notation~\ref{nvnot}, qui n'apparaît pas dans
\cite{artsmf}.

Si $X$ est un objet fini de $\F_{\Gr}$, le plus petit
$(n,k)$ tel que $X\in {\rm Ob}\,\F_{\Gr}[n,k]$ sera appelé {\em
  bidegré}\,\footnote{On prendra garde que la deuxième composante du
  bidegré d'un foncteur fini n'est pas toujours son degré !} de $X$ ;
on le note~${\rm bideg}\,X$. Il s'agit de la notion de taille d'un foncteur en grassmanniennes adaptés à cet article.

\paragraph*{Les foncteurs induits et le théorème de simplicité
  généralisé.}
Powell a introduit dans \cite{GP4} une filtration fondamentale du
foncteur décalage $\Delta_\FF : \F\to\F$ par des foncteurs noté
$\tilde{\nabla}_n$. Les foncteurs que nous notons ainsi sont en fait
les {\em duaux} de ceux introduits par Powell ; on a alors une suite
d'endofoncteurs
$$\Delta_\FF\simeq id\oplus\Delta\simeq\tilde{\nabla}_0\supset\tilde{\nabla}_1\simeq\Delta\supset\tilde{\nabla}_2\supset\dots\supset\tilde{\nabla}_n\supset\cdots$$
Nous renvoyons à \cite{GP4} pour les propriétés de ces foncteurs (qui
sont traduites dans notre contexte dual dans \cite{artsmf}, section~1.5). Rappelons
notamment que pour $n\geq
2$, le foncteur $\tilde{\nabla}_n$ n'est pas exact, mais il est
additif et préserve
toujours épimorphismes et monomorphismes ; par ailleurs, le foncteur
$\tilde{\nabla}_n$ annihile un foncteur simple $S_\lambda$ si et
seulement si la longueur $l(\lambda)$ de la partition régulière $\lambda$ est
strictement inférieure à~$n$. 

On peut montrer que la sous-catégorie pleine, notée $\N
il_{\tilde{\nabla}_n}$, des foncteurs $\tilde{\nabla}_n$-nilpotents de
$\F$ est épaisse (\cite{GP4}, §\,4.2).

\begin{theo}[Théorème de simplicité généralisé, \cite{artsmf}]\label{tsg} Pour tout  $n\in\mathbb{N}$, le foncteur $\omega_{n} : \F_{\Gr,n}\to\F$ induit une équivalence entre $\F_{\Gr,n}^f$ et une sous-catégorie épaisse
de $\N il_{\tilde{\nabla}_{n+1}}/\N il_{\tilde{\nabla}_n}$.
\end{theo}

Ce théorème (établi dans la section~16.2 de \cite{artsmf}) donne une indication forte en faveur des conjectures
équivalentes suivantes, discutées dans le chapitre~12 de \cite{artsmf}.

On rappelle que les notations $\mathbf{NT}_{n}(\F)$ (foncteurs noethériens de type $n$) et $\K_{n}(\F)$ (filtration de Krull de $\F$) qui apparaissent ci-après sont introduites au début de l'appendice.

\begin{conj}[Conjecture artinienne extrêmement forte, première version]\label{caef1} Pour tout  $n\in\mathbb{N}$, le foncteur $\omega_{n} : \F_{\Gr,n}\to\F$ induit une équivalence entre
$\F_{\Gr,n}^f$ et $\mathbf{NT}_{n}(\F)/\mathbf{NT}_{n-1}(\F)$.
\end{conj}

\begin{conj}[Conjecture artinienne extrêmement forte, deuxième version]\label{caef2} Pour tout  $n\in\mathbb{N}$, le foncteur $\omega_{n} : \F_{\Gr,n}\to\F$ induit une équivalence entre
$\F_{\Gr,n}^{lf}$ et $\K_{n}(\F)/\K_{n-1}(\F)$.
\end{conj}

Le principe intuitif de la conjecture artinienne extrêmement forte
consiste à ramener tous les foncteurs de type fini de la catégorie $\F$ à
l'image par le foncteur $\omega$ de foncteurs finis de la catégorie
$\F_\Gr$. Mais tout objet de $\F^{tf}$ (i.e. de type fini --- cf. appendice) n'appartient pas à l'image
essentielle de la restriction de $\omega$ aux objets finis. Pour
donner une dernière forme de la conjecture qui ne fasse pas apparaître
de catégorie quotient, nous sommes ainsi amenés à introduire la définition suivante.

\begin{defi}\label{df-sep}\begin{itemize}\item Disons qu'une sous-catégorie pleine $\C$ d'une catégorie abélienne $\A$ est
{\em semi-épaisse} si elle contient $0$ et que pour toute suite exacte
courte $0\to A\to B\to C\to 0$ de $\A$, si deux des trois objets de
$A$, $B$ et $C$ sont dans $\C$, alors il en est de même pour le troisième.
\item Pour tout entier $n$, on note $\F^{\omega - cons(n)}$ la
  sous-catégorie semi-épaisse de $\F$ engendrée par les $\omega(X)$, où
  $X\in {\rm Ob}\,\F_\Gr^f$ est de bidegré~$(k,j)$ avec $k\leq n$. Les
  objets de $\F^{\omega - cons(n)}$ sont appelés {\em foncteurs
    induits de hauteur au plus~$n$} de~$\F$.
\end{itemize}\end{defi}

On rappelle qu'une sous-catégorie épaisse est une sous-catégorie semi-épaisse qui est stable par sous-objets.

La catégorie $\F^{\omega - cons(n)}$ est constituée des foncteurs qui sont constructibles à partir des foncteurs $\omega(X)$ du type indiqué dans l'énoncé en ce sens qu'on peut les obtenir en un nombre fini d'étapes du type extension, noyau d'épimorphisme, conoyau de monomorphisme à partir de ces derniers. Cela justifie la notation.

\begin{conj}[Conjecture artinienne extrêmement forte, troisième version]\label{caef3} Pour tout  $n\in\mathbb{N}$, la sous-catégorie $\F^{\omega - cons(n)}$ de $\F$ est
épaisse.
\end{conj}

Cette conjecture s'aborde par récurrence sur l'entier $n$, ce qui amène à introduire l'hypothèse suivante :

\bigskip

 $H^\omega(n)$ \qquad {\em Pour tout entier $i\leq n$, la sous-catégorie $\F^{\omega - cons (i)}$ de $\F$
  est épaisse.}
  
  \bigskip
  
La première assertion de la proposition suivante est un corollaire
direct du théorème de simplicité généralisé ; pour les autres points,
voir \cite{artsmf}, section~12.2.

\begin{pr}\label{cortsg} Soit $n\in\mathbb{N}$. \begin{enumerate}\item Si l'hypothèse $H^\omega(n-1)$ est
satisfaite et que tout foncteur de type fini
$\tilde{\nabla}_n$-nilpotent appartient à $\F^{\omega-cons(n-1)}$, alors
l'hypothèse $H^\omega(n)$ est satisfaite .
\item Si l'hypothèse $H^\omega(n)$ est satisfaite, alors le foncteur $\omega_k : \F_{\Gr,k}\to\F$ induit une équivalence de
catégories $\F_{\Gr,k}^f\to\F^{\omega - cons(k)}/\F^{\omega - cons(k-1)}$
pour $k\leq n$.
\item Si l'hypothèse $H^\omega(n)$ est satisfaite, alors $P^{\otimes i}\otimes F$ est noethérien de type $i$ pour tout foncteur
fini $F$ de $\F$ et tout entier naturel $i\leq n$.
\end{enumerate}
\end{pr}

\begin{rem}\label{rqbideg} Supposons l'hypothèse $H^\omega(n)$ satisfaite. La deuxième assertion
de la proposition~\ref{cortsg} permet de transporter la notion de
bidegré $\F^{\omega - cons(n)}$ : pour $F\in {\rm Ob}\,\F^{\omega - cons(n)}$, on définit 
le {\em bidegré} de $F$, noté $||F||$, comme le couple $(k,j)$, où $k$
est le plus petit entier tel que $F\in {\rm Ob}\,\F^{\omega - cons(k)}$,
et $j$ le degré du foncteur fini $X$ de $\F_{\Gr,k}$ tel que
$F\simeq\omega_k(X)$ dans $\F^{\omega - cons(k)}/\F^{\omega -
  cons(k-1)}$. Ainsi, pour tout $Y\in {\rm Ob}\,\F_\Gr^f$ de bidegré
$(k,j)$ avec $k\leq n$, le foncteur $\omega(X)$ de $\F$ est de bidegré
$(k,j)$.
\end{rem}

\paragraph*{Stratégie d'approche de la conjecture artinienne
  extrêmement forte.}  Le fait que $H^\omega(1)$ est vérifiée est une
conséquence directe du théorème de simplicité généralisé et de ce que
les objets finis de $\F$ sont polynomiaux ; dans \cite{artsmf} (section~16.4),
  nous avons observé que les résultats de Powell sur les foncteurs annulés par $\tilde{\nabla}_{2}$ (cf. \cite{GP5}) impliquent $H^\omega(2)$ (résultat que nous retrouverons
  d'ailleurs par une méthode différente dans le présent
  article). Notre but consiste à établir $H^\omega(3)$. Dans la
  section suivante, nous ramènerons la vérification de cette hypothèse
  à celle d'énoncés relatifs à des foncteurs de division, grâce à la proposition~\ref{pr-cadiv}.

\section{Étude des foncteurs de division}\label{s-div}

Si $F$ est un foncteur de $\F$ à valeurs de dimension finie, alors
l'endofoncteur $-\otimes F$ de $\F$ commute aux limites ; pour des
raisons formelles, on en déduit qu'il possède un adjoint à droite,
appelé foncteur de {\em division par $F$} et noté $(- : F)$ (cf. \cite{artsmf}, appendice C3). En fait,
$(-:-)$ définit même un bifoncteur $\F\times\F^{op}\to\F$ exact à
droite en chaque variable.
Pour des
raisons typographiques, le foncteur $(- : F)$  sera aussi noté $\D[F]$,
et ses foncteurs dérivés à gauche $\D^i[F]$. Une observation évidente,
omniprésente dans la suite, réside dans la commutation naturelle de
deux foncteurs de division : $\D[F]\circ\D[G]\simeq\D[F\otimes
G]\simeq\D[G]\circ\D[F]$ ; en particulier, la $i$-ième itérée
$\D[F]^{\circ i}$ (notation à ne pas confondre avec $\D^i[F]$ !) de
$\D[F]$ est isomorphe à $\D[F^{\otimes i}]$. La division par un
foncteur injectif (à valeurs de dimension finie) de $\F$ est un
foncteur {\em exact} ; la division par un foncteur injectif de co-type
fini\,\footnote{Cf. début de l'appendice pour cette notion.} commute également aux limites.

On peut donner une description explicite des foncteurs de division à
partir du lemme de Yoneda (cf. \cite{art1} et \cite{artsmf}) : on a $(P_V
: F)\simeq F(V)^*\otimes P_V$ (ce qui montre d'ailleurs que les
foncteurs de division préservent $\F^{tf}$) et
\begin{equation}\label{div-expl}(A:B)(V)^*\simeq {\rm hom}_\F(\Delta_V(A),B).\end{equation}
Les foncteurs $\D^i[B]$ se décrivent de façon similaire en remplaçant
hom par Ext$^i$.

\begin{rem}[fondamentale] {\em Les foncteurs $\tilde{\nabla}_n$ peuvent se définir à partir de
foncteurs de division} : il existe des foncteurs $A_n$ et $B_n$ à
valeurs de dimension finie et un morphisme $f_n : B_n\to A_n$ de $\F$
tels que $\tilde{\nabla}_n=im\,(\D[A_n]\xrightarrow{\D[f_n]}\D[B_n])$.
\end{rem}

Il existe de même, dans la catégorie $\F_\Gr$, des foncteurs de
division $(-:X)_\Gr$, où $X\in {\rm Ob}\,\F_\Gr$ prend des valeurs de
dimension finie. On adoptera aussi la notation $\D_\Gr[X]$, et
$\D^i_\Gr[X]$ pour les foncteurs dérivés. Dans $\F_{\Gr,n}$ et $\F_{\Gr,\leq n}$, on emploie les notations analogues $\D^i_{\Gr,n}[X]$ et $\D^i_{\Gr,\leq n}[X]$.

Une observation élémentaire mais fondamentale est que, dans $\F$ comme dans $\F_{\Gr}$ ou $\F_{\Gr,n}$, {\em le foncteur différence est un foncteur de division} (le cas de $\F$ est traité dans l'appendice de
\cite{GP4}, celui des catégories de foncteurs en grassmanniennes dans
la section~5.3 de \cite{artsmf}) ; par conséquent,
tous les foncteurs de division commutent, à isomorphisme naturel près,
au foncteur différence. Comme le foncteur différence est exact et préserve les
foncteurs projectifs, on en déduit le résultat suivant.

\begin{pr}\label{prcomd} Les dérivés des foncteurs de
division commutent, à isomorphisme naturel près, au foncteur différence.
\end{pr}

\subsection{Lien entre les foncteurs $\D[S_{<n>}]$, $\D[L(n)]$ et $\tilde{\nabla}_n$}

L'objectif de ce paragraphe consiste à montrer que les foncteurs $\D[S_{<n>}]$, $\D[L(n)]$ (les notations $<n>$, $S_{<n>}$ et $L(n)$ sont introduites dans l'appendice, remarque~\ref{rqstein}) et $\tilde{\nabla}_n$ annihilent les mêmes foncteurs simples de $\F$, ce qui permettra notamment
de remplacer la condition de $\tilde{\nabla}_n$-nilpotence à manier pour vérifier l'hypothèse $H^\omega(n)$ par une condition de $\D[L(n)]$-nilpotence.

Tous les cas de non-annulation par division considérés dans cet article s'appuient sur les deux lemmes suivants, dont la démonstration
repose sur des considérations explicites sur les foncteurs de Weyl et les foncteurs simples de~$\F$.

\begin{lm}\label{lm1} Soient $\lambda$ et $\mu$ deux partitions telles que $\lambda+\mu$ soit régulière. Il existe un épimorphisme $W_\lambda\otimes
  W_\mu\twoheadrightarrow W_{\lambda+\mu}$, où $\lambda+\mu$ désigne
  la partition $(\lambda_1+\mu_1,\dots,\lambda_i+\mu_i,\dots)$.
\end{lm}

\begin{proof} Cette propriété s'obtient aussitôt à l'aide des éléments
  semi-standard des foncteurs de Weyl, qui engendrent $W_{\lambda+\mu}$ par l'hypothèse de régularité de $\lambda+\mu$ (cf. remarque~\ref{rq-gst}). L'épimorphisme recherché
  est la restriction de la flèche
  $\Lambda^\lambda\otimes\Lambda^\mu\twoheadrightarrow\Lambda^{\lambda+\mu}$ obtenue par produit tensoriel des morphismes canoniques (produits) $\Lambda^{\lambda_i}\otimes\Lambda^{\mu_i}\twoheadrightarrow\Lambda^{\lambda_i+\mu_i}$.
\end{proof}
                   
\begin{lm}\label{lm2} Soit $\lambda$ une partition régulière de
  longueur au moins $n$. Il existe une partition $\alpha$ et un
  épimorphisme $W_\alpha\otimes S_{<n>}\twoheadrightarrow S_\lambda$.
\end{lm}

\begin{proof} On rappelle d'abord que $S_{<n>}\simeq W_{<n>}$
  (cf. remarque~\ref{rqstein}). Comme la partition $\lambda$ est régulière,
  il existe une partition $\alpha$  telle que $\alpha+<n>=\lambda$. On obtient donc par le lemme~\ref{lm1} un épimorphisme $W_\alpha\otimes S_{<n>}\twoheadrightarrow W_{\alpha+<n>}=W_\lambda\twoheadrightarrow S_\lambda$.
\end{proof}

La proposition suivante justifie, à l'aune du théorème de simplicité
généralisé, l'intérêt pour les foncteurs de division par $S_{<n>}$ et
$L(n)$. L'avantage spécifique du premier
est d'être~\guillemotleft~plus petit~\guillemotright, celui du second réside dans son exactitude.

\begin{pr}\label{pr-divst} Soient $\lambda$ une partition régulière et
  $n\in\mathbb{N}$. Les
  assertions suivantes sont équivalentes : 
\begin{enumerate}\item $l(\lambda)<n$ (on rappelle que $l(\lambda)$ désigne la longueur de $\lambda$ --- cf. appendice) ;
\item $(S_\lambda:L(n))=0$ ;
\item $(S_\lambda:S_{<n>})=0$ ;
\item $\tilde{\nabla}_n(S_\lambda)=0$.
\end{enumerate}
\end{pr}

\begin{proof} L'équivalence entre les assertions 1 et 4 est déjà
  connue (cf. \cite{GP4}, théorème~1).

Grâce à l'isomorphisme (\ref{div-expl}), si l'assertion 2
est fausse, il existe $V\in {\rm Ob}\,\E^f$ et un morphisme non nul
$\Delta_V(S_\lambda)\to L(n)$, de sorte que $\Delta_V(S_\lambda)$ a un facteur
de composition $S_{<n>}$. En appliquant le
foncteur $\tilde{\nabla}_n$, qui préserve  épimorphismes et monomorphismes, commute
aux foncteurs de décalage, et n'annihile pas $S_{<n>}$, on en déduit que
$\tilde{\nabla}_n(S_\lambda)$ est non nul. Ainsi $4\Rightarrow 2$. 

L'implication $2\Rightarrow 3$ est évidente puisque l'inclusion
$S_{<n>}\hookrightarrow L(n)$ induit une surjection $\D[L(n)]\twoheadrightarrow\D[S_{<n>}]$.

Montrons que 3 entraîne 1. Si $\lambda$ est de longueur $\geq n$,
alors on dispose par le lemme~\ref{lm2} d'un épimorphisme
$W_\alpha\otimes S_{<n>}\twoheadrightarrow S_\lambda$, d'où par
dualité $S_\lambda\hookrightarrow S_{<n>}\otimes DW_\alpha$ (les
objets simples de $\F$ étant auto-duaux), ce qui
montre que $(S_\lambda : S_{<n>})$ est non nul.
\end{proof}

\begin{rem}\label{rqdiv} On a $\tilde{\nabla}_n(F)=0\;\Rightarrow\;
  (F:L(n))=0\;\Rightarrow\; (F:S_{<n>})=0$ pour tout  $F\in {\rm Ob}\,\F$. 
\end{rem}

Nous déduisons maintenant de la proposition~\ref{pr-divst} la comparaison entre
foncteurs $\tilde{\nabla}_n$-nilpotents et $\D[L(n)]$-nilpotents. 

\begin{cor}\label{cr-andiv} Soient $F\in {\rm Ob}\,\F$ et $n\in\mathbb{N}$.
\begin{enumerate}
\item Les assertions suivantes sont équivalentes :
\begin{enumerate}
\item $(F:L(n))=0$ ;
\item si un foncteur simple $S_{\lambda}$ est facteur de composition de $F$, alors la partition régulière
$\lambda$ est de longueur strictement inférieure à~$n$ ;
\item tous les facteurs de composition de $F$ sont annihilés par le foncteur $\tilde{\nabla}_{n}$.
\end{enumerate}
\item Si le foncteur $F$ est annulé par une itérée $\tilde{\nabla}_{n}^{\circ i}$ de $\tilde{\nabla}_{n}$, il est annulé par l'itérée $\D[L(n)]^{\circ i}$ de $\D[L(n)]$.
\end{enumerate}
\end{cor}

\begin{proof}\begin{enumerate}\item  Comme le foncteur $\D[L(n)]$ est
    exact et commute
    aux limites et colimites, la condition $\D[L(n)](F)=0$ équivaut à
    $\D[L(n)](S_\lambda)=0$ pour tout facteur de composition
    $S_\lambda$ de $F$ (écrire ce foncteur comme colimite de foncteurs
    de type fini, qui sont limites de foncteurs finis). L'équivalence
    des trois conditions résulte alors de la proposition~\ref{pr-divst}.
\item Si $\tilde{\nabla}_n(F)=0$, alors
  $\tilde{\nabla}_{n}(S_\lambda)=0$ pour tout facteur de composition
  $S_\lambda$ de $F$, car le foncteur $\tilde{\nabla}_n$ préserve les
  monomorphismes et les épimorphismes. Ce qui précède montre alors que
  $\tilde{\nabla}_n(F)=0$ entraîne $\D[L(n)](F)=0$, ce qui traite le
  cas $i=1$. Le cas général s'en déduit par récurrence, parce que les
  foncteurs $\tilde{\nabla}_n$ et $\D[L(n)]$ commutent, par le lemme~\ref{lmcom} ci-après.
\end{enumerate}
\end{proof}

\begin{lm}\label{lmcom} Soient $n$ et $m$ deux entiers naturels. Les
  foncteurs $\tilde{\nabla}_n$ et $\D[L(m)]$ commutent à isomorphisme près.
\end{lm}

\begin{proof} \'Ecrivons $\tilde{\nabla}_n=im\,(\D[A]\xrightarrow{\D[f]}\D[B])$ où
  $B\xrightarrow{f}A$ est un morphisme convenable de $\F$ (avec $A$ et
  $B$ à valeurs de dimension finie). Il existe un diagramme
  commutatif
$$\xymatrix{\D[A]\circ\D[L(m)]\ar[rr]^-{\D[f]\D[L(m)]}\ar[d]_\simeq &
  & \D[B]\circ\D[L(m)]\ar[d]^\simeq \\
\D[L(m)]\circ\D[A]\ar[rr]^-{\D[L(m)]\D[f]} & & \D[L(m)]\circ\D[B].
}$$
L'image de la flèche supérieure est $\tilde{\nabla}_n\circ\D[L(m)]$,
et celle de la flèche inférieure s'identifie à
$\D[L(m)]\circ\tilde{\nabla}_n$ par exactitude du foncteur
$\D[L(m)]$. Cela démontre le lemme.
\end{proof}

\subsection{Foncteurs de division dans $\F_{\Gr}$}\label{par-dgr}

Nous commençons par établir un lien formel entre les foncteurs de division dans $\F$ et dans $\F_{\Gr}$. Il s'agit d'étudier l'effet de foncteurs de division sur des foncteurs
de type fini de $\F$ à partir de l'effet de foncteurs de division sur des foncteurs {\em finis} de $\F_{\Gr}$, via le foncteur $\omega$.

Le contenu de ce paragraphe recoupe partiellement les considérations
de \cite{artsmf} (chapitre~9) sur les foncteurs de division.

\begin{pr}\label{pr1-dom} Il existe un isomorphisme naturel gradué
$$\D^*[F]\circ\omega\simeq \omega\circ\D^*_\Gr[\iota(F)]$$
de foncteurs $\F_\Gr\to\F$, pour $F\in {\rm Ob}\,\F$ à valeurs de
dimension finie.
\end{pr}

\begin{proof} Le cas du degré $0$ est traité dans la section~9.2 de \cite{artsmf} (il
  s'agit d'une conséquence formelle de
  l'adjonction entre les foncteurs $\iota$ et $\omega$) ; le cas général s'en déduit aussitôt
  parce que le foncteur $\omega$ est exact et préserve les objets
  projectifs, son adjoint à droite $\iota$ étant exact.
\end{proof}

Les deux propositions suivantes permettent de passer de la catégorie $\F_{\Gr}$ aux catégories $\F_{\Gr,n}$.

\begin{pr}\label{prdrp} Soient $X$ un objet de $\F_\Gr$  à valeurs de
  dimension finie et $n\in\mathbb{N}$.
Il existe un isomorphisme naturel gradué
$\D_\Gr^*[X]\circ\Prl_{\leq n}\simeq\Prl_{\leq n}\circ\D^*_{\Gr,\leq n}[\R_{\leq n}(X)]$
de foncteurs $\F_{\Gr,\leq n}\to\F_\Gr$.
\end{pr}

\begin{proof} Le foncteur de prolongement par zéro $\Prl_{\leq n} : \F_{\Gr,\leq n}\to\F_\Gr$
  est adjoint à gauche au foncteur de restriction $\R_{\leq n} :
  \F_\Gr\to\F_{\Gr, \leq n}$  (cf. \cite{artsmf}, section~5.1), qui commute au produit tensoriel. Par conséquent, les foncteurs
  $\D_\Gr^*[X]\circ\Prl_{\leq n}$ et $\Prl_{\leq n}\circ\D^*_{\Gr,\leq
    n}[\R_{\leq n}(X)]$ sont isomorphes : ils sont tous deux adjoints
  à gauche au foncteur composé $\R_{\leq n}\circ (-\otimes X)$, ce qui
  fournit l'isomorphisme recherché
  en degré $0$. Le cas général s'en déduit, comme dans la
  démonstration de la proposition~\ref{pr1-dom}, par exactitude des foncteurs
  $\Prl_{\leq n}$ et $\R_{\leq n}$.
\end{proof}

Dans la proposition suivante, on note, par abus, $\R_n :
\F_{\Gr,\leq n}\to\F_{\Gr,n}$ et $\Prl_n :
\F_{\Gr,n}\to\F_{\Gr,\leq n}$ les foncteurs de restriction et de
prolongement par zéro respectivement, analogues à ceux considérés dans
la section~\ref{s1}.

\begin{pr}\label{prdrp2} Soient $n\in\mathbb{N}$ et $X$ un objet de $\F_{\Gr,\leq n}$ à valeurs de dimension finie.
 Il existe un isomorphisme naturel gradué
$\R_n\circ\D^*_{\Gr,\leq n}[X]\simeq\D^*_{\Gr,n}[\R_n(X)]\circ\R_n$ de
foncteurs $\F_\Gr\to\F_{\Gr,n}$.
\end{pr}

\begin{proof} Cette fois-ci, le foncteur de prolongement par zéro $\Prl_n : \F_{\Gr, n}\to\F_{\Gr,\leq n}$
  est adjoint à {\em droite} au foncteur de restriction $\R_{\leq n} :
  \F_{\Gr, \leq n}\to\F_{\Gr,n}$ (cf. \cite{artsmf}, §\,5.1). On conclut par un
  argument formel similaire à celui invoqué pour la proposition~\ref{prdrp}.
\end{proof}

Nous examinons maintenant le comportement des foncteurs de division dans $\F_{\Gr,n}$ sur les foncteurs {\em pseudo-constants} (ces foncteurs sont définis en~\ref{rqf} ; les notations $\rho_n$ et $\epsilon_n$ qui apparaissent dans l'énoncé ci-dessous sont celles introduites en~\ref{dfpc}).

\begin{pr}\label{pr-and} Soient $n\in\mathbb{N}$, $X$ un objet de
$\F_{\Gr,n}$ à valeurs de dimension finie et $M$ un $GL_n$-module. Il
existe un isomorphisme naturel
$$\D^i_{\Gr,n}[X](\rho_n(M))\simeq\left\{\begin{array}{cc}\rho_n(M\otimes\varepsilon_n(X)^*)
& \text{si } i=0,\\
0 & \text{si }i>0.
\end{array}\right.$$
\end{pr}

\begin{proof} La catégorie monoïdale symétrique  $GL_n - \mathbf{mod}$
  (le produit tensoriel étant pris sur $\FF$) possède également des
  foncteurs de division par des modules finis, qui sont donnés par
  $(M:N)_{GL_n}=M\otimes N^*$ (cf. \cite{CR}, §\,10.D, par exemple). En
  particulier, ces foncteurs sont {\em exacts}.

Maintenant, la proposition s'obtient formellement (cf. les
démonstrations précédentes) à partir des trois
observations suivantes :
\begin{itemize}\item le foncteur $\rho_n$ est adjoint à gauche à
  $\varepsilon_n$ (cf. proposition \ref{prpc}) ;
\item le foncteur $\varepsilon_n$ commute au produit tensoriel ;
\item les foncteurs $\rho_n$ et $\varepsilon_n$ sont exacts.
\end{itemize}
\end{proof}

La proposition suivante constitue le résultat principal sur la
catégorie $\F_\Gr$ de ce
paragraphe.

\begin{pr}\label{pr-tdgr} Soit $X$ un objet de $\F_\Gr$ à valeurs de
dimension finie.
\begin{enumerate}\item Les foncteurs $\D^i_\Gr[X]$ préservent les
  objets finis.
\item Pour tous $i\in\mathbb{N}$ et $Y\in {\rm Ob}\,\F_\Gr^f$, on a
  ${\rm bideg}\,\D^i_\Gr[X](Y)\leq {\rm bideg}\,Y$. De plus, si $Y$
  est non nul, l'égalité advient si et seulement si les deux
  conditions sont satisfaites :
\begin{enumerate}\item $i=0$ ;
\item $X(\FF^n,\FF^n)\neq 0$, où $n$ est la première composante de ${\rm bideg}\,Y$.
\end{enumerate}
\end{enumerate}
\end{pr}

\begin{proof} Nous avons observé que
  les foncteurs de division et leurs dérivés commutent au foncteur
  différence (proposition~\ref{prcomd}). Par conséquent, les foncteurs $\D^i_\Gr[X]$ préservent
  les foncteurs polynomiaux et diminuent (au sens large) leur degré.

D'autre part, pour $Y\in {\rm Ob}\,\F_\Gr^f$, les $\D^i_\Gr[X](Y)$
sont à valeurs de dimension finie. C'est une conséquence facile du
lemme de Yoneda et de ce que les objets finis de $\F_\Gr$ possèdent
des résolutions projectives de type fini (cf. \cite{artsmf}, §\,5.5).

Si $Y\in {\rm Ob}\,\F_\Gr^f$ est non nul et de bidegré $(n,k)$, il
existe $Z\in {\rm Ob}\,\F_{\Gr,\leq n}$ tel que $Y\simeq\Prl_{\leq
  n}(Z)$, et que $\R_n(Z)=\R_n(Y)$ est de degré $k$. Par la proposition~\ref{prdrp}, on a
$\D_\Gr^i[X](Y)\simeq\Prl_{\leq n}\D_{\Gr,\leq n}^i[\R_{\leq n}(X)](Z)$.

Cela prouve déjà, compte-tenu du début de la démonstration et de la
proposition~\ref{pol-fgr}, que les foncteurs $\D^i_\Gr[X]$ préservent les
  objets finis. On en déduit également que le bidegré de
$\D_\Gr^i(Y)$ est $(n,\deg \R_n\D_{\Gr,\leq n}^i[\R_{\leq n}(X)](Z))$ si
$\R_n\D_{\Gr,\leq n}^i[\R_{\leq n}(X)](Z)\neq 0$, et strictement inférieur à $(n,k)$
si $\R_n\D_{\Gr,\leq n}^i[\R_{\leq n}(X)](Z)=0$.

Par la proposition~\ref{prdrp2}, $\R_n\D_{\Gr,\leq n}^i[\R_{\leq
  n}(X)](Z)\simeq\D^i_{\Gr,n}[\R_n(X)](\R_n(Y))$. Si $\R_n(Y)$ est
pseudo-constant (i.e. $k=0$), la proposition~\ref{pr-and} donne la
conclusion souhaitée. Le cas général s'en déduit par récurrence sur
$k$, puisque les foncteurs $\D^i_{\Gr,n}[X]$ commutent au foncteur différence.
\end{proof}

\begin{cor}\label{cr-doa} Soient $n\in\mathbb{N}$ et $F\in {\rm
    Ob}\,\F$ à valeurs de dimension finie tel que $F(\FF^n)=0$. Alors tout
  foncteur induit de hauteur au plus $n$ est $\D[F]$-nilpotent.
\end{cor}

\begin{proof} La proposition~\ref{pr-tdgr} montre que si $X$ est un
  foncteur fini de $\F_\Gr$ dont la première composante du bidegré est
  au plus~$n$, alors il existe un entier $k$ pour lequel
  $$\D^{i_1}_\Gr[\iota(F)]\D^{i_2}_\Gr[\iota(F)]\dots\D^{i_k}_\Gr[\iota(F)](X)=0$$ pour toute suite d'entiers naturels $(i_1,i_2,\dots,i_k)$, puisque $\iota(F)(\FF^n,\FF^n)=F(\FF^n)=0$ par hypothèse.

On en déduit, par la proposition~\ref{pr1-dom}, que le foncteur
$A=\omega(X)$ est tel que
 $$\D^{i_1}[F]\D^{i_2}[F]\dots\D^{i_k}[F](A)=0$$ pour toute suite
 d'entiers naturels $(i_1,i_2,\dots,i_k)$.

Le corollaire découle alors de ce que la sous-catégorie pleine des
foncteurs $A$ de $\F$ pour
lesquels existe un tel $k$ est semi-épaisse (définition~\ref{df-sep}).
\end{proof}


%

Rappelons que la notion de bidegré noté $||.||$ dans la catégorie $\F$ a été introduite dans la remarque~\ref{rqbideg}.

\begin{pr}\label{pr-tdgrf}Soient $n\in\mathbb{N}$ tel que l'hypothèse $H^\omega(n)$ est
vérifiée, $F$ un foncteur de $\F$ à valeurs de dimension finie et $G\in {\rm Ob}\,\F^{\omega-cons(n)}$.
\begin{enumerate}\item Pour tout $i\in\mathbb{N}$, on a $\D^i[F](G)\in {\rm Ob}\,\F^{\omega-cons(n)}$.
\item On a $||\D^i[F](G)||\leq ||G||$. Lorsque $G$ est non nul,
  l'égalité a lieu si et seulement si les deux conditions suivantes
  sont satisfaites :
\begin{enumerate}\item $i=0$ ;
\item $F(\FF^j)\neq 0$, où $j$ est la première composante de $||G||$.
\end{enumerate}
\end{enumerate}
\end{pr}

\begin{proof} On raisonne par récurrence sur $||G||$, noté $(j,k)$ : on suppose le
  résultat connu pour les bidegrés strictement inférieurs à $(j,k)$. Soit $X\in {\rm
    Ob}\,\F^f_{\Gr,j}$ un foncteur de degré $k$ tel que
  $G\simeq\omega_j(X)$ dans la catégorie $\F^{\omega-cons(j)}/\F^{\omega-cons(j-1)}$. Il
  existe donc une suite exacte
$$0\to N\to G\xrightarrow{f}\omega_j(X)\to C\to 0$$
dans laquelle les foncteurs $N$ et $C$ appartiennent à
  $\F^{\omega -cons (j-1)}$. Par conséquent, leurs bidegrés sont
  strictement inférieurs à $(j,k)$, et l'hypothèse de récurrence
donne en particulier $||\D^i[F](N)||<(j,k)$ et
$||\D^i[F](C)||<(j,k)$ pour tout $i$. Par ailleurs, la
proposition~\ref{pr1-dom} donne
$\D^i[F](\omega_j(X))\simeq\omega(\D^i_\Gr[\iota(F)](Y))$, où
$Y=\Prl_n(X)\in {\rm Ob}\,\F_\Gr^f$ est de bidegré $(k,j)$. La
proposition~\ref{pr-tdgr} décrit donc le bidegré de $\D^i[F](\omega_j(X))$.

Notons $H=im\,f$ : les suites exactes longues pour $\D^*[F]$ associées aux suites exactes
courtes $0\to H\to\omega_j(X)\to C\to 0$ puis $0\to N\to G\to H\to 0$
donnent alors la conclusion, puisque dans une suite exacte $0\to A\to
B\to C\to 0$ de $\F^{\omega-cons(n)}$, on a
$||B||=\max(||A||,||C||)$.
\end{proof}

\subsection{Foncteurs de division et foncteurs induits}

L'hypothèse suivante, dans laquelle $A$ désigne un foncteur de $\F$ à valeurs de dimension finie et $i$ un entier, jouera un rôle fondamental dans cet article.

\bigskip

$H_{div}(A,i)$\qquad {\em Si $F\in {\rm Ob}\,\F^{tf}$ vérifie $(F:A)=0$, alors $F\in {\rm Ob}\,\F^{\omega - cons(i-1)}$.}

\bigskip

\begin{rem}\label{rqft}\begin{enumerate}
\item Cette hypothèse ne peut être vérifiée que si $A(\FF^{i})\neq 0$, puisque le foncteur projectif $P_{\FF^{i}}$ n'est pas dans $\F^{\omega - cons(i-1)}$ et que
$(P_{\FF^{i}}:A)\simeq A(\FF^{i})^*\otimes P_{\FF^{i}}$. Si la conjecture artinienne extrêmement forte est vraie, alors
réciproquement, pour tout foncteur $A$ à valeurs de dimension finie tel que $A(\FF^{i})\neq 0$, l'hypothèse $H_{div}(A,i)$ est satisfaite (cela découle de la proposition~\ref{pr-tdgrf}).
\item L'hypothèse $H_{div}(A,i)$ est surtout intéressante lorsque le
  foncteur $\D[A]$ est nilpotent sur $\F^{\omega - cons(i-1)}$,
  condition équivalente à $A(\FF^{i-1})=0$ (cf. le paragraphe précédent,
  notamment le corollaire~\ref{cr-doa}).
\item Nous ne manierons l'hypothèse $H_{div}(A,i)$ que lorsque l'hypothèse $H^{\omega}(i-1)$ est satisfaite, le but de la considération de $H_{div}(A,i)$ consistant alors
à établir $H^\omega(i)$.
\item\label{itf} Si $A$ est un sous-foncteur d'un produit tensoriel $B\otimes C$
  de foncteurs à valeurs de dimension
  finie $B\otimes C$, alors l'hypothèse $H_{div}(A,i)$ entraîne $H_{div}(B,i)$,
  puisque l'inclusion induit un épimorphisme $\D[C]\D[B]\simeq\D[B\otimes C]\twoheadrightarrow\D[A]$.
\item  Le cas où $A$ est un foncteur injectif, équivalent à
  l'exactitude de $\D[A]$, joue un rôle particulier (cf. lemme~\ref{lm-cld}). Nous aurons à manier
l'hypothèse $H_{div}(A,i)$ pour différents foncteurs $A$, notamment des foncteurs simples, mais la connaissance de la situation pour certains foncteurs injectifs sera nécessaire.
\item Les remarques précédentes illustrent l'intérêt de $H_{div}(A,i)$
  pour $A$ injectif tel que $A(\FF^{i-1})=0$ et $A(\FF^{i})\neq 0$. Le
  plus petit --- en ce sens que la hauteur de
    son dual, comme
    foncteur induit, est minimale --- foncteur $A$ de ce
type est le foncteur $L(i)$ ; c'est pourquoi l'hypothèse $H_{div}(L(i),i)$ occupera une place centrale dans cet article.
\end{enumerate}\end{rem}

\begin{nota} Nous désignerons par
$H_{div}(i)$ l'hypothèse $H_{div}(L(i),i)$.
\end{nota}

\begin{lm}\label{lm-cld} Sous les hypothèses $H_{div}(i)$ et $H^\omega(i-1)$, tout
  foncteur de type fini $\D[L(i)]$-nilpotent est induit de
  hauteur au plus~$i-1$. 
\end{lm}

\begin{proof} Soient $F\in {\rm Ob}\,\F^{tf}$, $f: F\to\D[L(i)](F)\otimes L(i)$ l'unité de l'adjonction
  et $N=ker\,f$. On a $\D[L(i)](N)=0$, car $\D[L(i)]$ est
  exact et $\D[L(i)](f)$ injectif. Comme $F$ est de type fini et
  $L(i)$ localement fini, il existe un
  sous-foncteur fini $A$ de $L(i)$ tel que
  $im\,f\subset\D[L(i)](F)\otimes A$. 

Si  $\D[L(i)](F)$ est
  induit de hauteur au plus $i-1$, la stabilité de $\F^{\omega -
    cons (i-1)}$ par $-\otimes A$ (cf. \cite{artsmf}, §\,12.1) et l'hypothèse  $H^\omega(i-1)$
  assurent que $im\,f$ appartient à $\F^{\omega -
    cons (i-1)}$. En particulier, comme tout foncteur induit est
  de présentation finie (cf. ibid.), $N$ est de type
  fini, et l'hypothèse $H_{div}(i)$ assure
  que $N$ est induit de hauteur au plus $i-1$, donc $F$
  aussi. Une récurrence immédiate fournit alors la conclusion.
\end{proof}

Les progrès sur la conjecture artinienne présentés dans cet article reposent sur la proposition suivante.

\begin{pr}\label{pr-cadiv} Supposons que les hypothèses $H_{div}(i)$ et $H^\omega(i-1)$ soient vérifiées. Alors il en est de même pour  $H^\omega(i)$.
\end{pr}

\begin{proof} Ce résultat découle du corollaire~\ref{cr-andiv}, du
  lemme~\ref{lm-cld} et de la proposition~\ref{cortsg}.
\end{proof}

La dernière section de ce travail démontrera la validité de
l'hypothèse $H_{div}(i)$ pour~$i\leq 3$.

L'énoncé suivant, qui s'appuie sur le lemme~\ref{lm-cld} et le
paragraphe~\ref{par-dgr}, permettra de mener à bien le pas de la
récurrence intervenant dans l'argument principal de cet article.

\begin{pr}\label{pr-omd} Soient $i$, $j$ et $n$ des entiers tels que $i-1\leq j\leq n$ et que l'hypothèse  $H^\omega(n)$ soit vérifiée, et $A\in {\rm Ob}\,\F$ un
foncteur à valeurs de dimension finie. On suppose que les hypothèses $H^\omega(n)$, $H_{div}(i)$ et $H_{div}(A,i)$ sont vérifiées.

Si $F\in {\rm Ob}\,\F^{tf}$ est tel que
$(F:A)\in {\rm Ob}\,\F^{\omega - cons(j)}$, alors $F\in {\rm Ob}\,\F^{\omega - cons(j)}$.
\end{pr}

\begin{proof} On procède par deux récurrences imbriquées, dont la
  première porte sur~$j$.

On commence par traiter le cas $j=i-1$. Comme tout
  foncteur  induit de hauteur au plus $i-1$ est $\D[L(i)]$-nilpotent par
  le corollaire~\ref{cr-doa}, il
  existe $k$ tel que $\D[L(i)]^{\circ k}\D[A](F)=0$. On a donc
  $\D[A]\D[L(i)]^{\circ k}(F)=0$, puis $\D[L(i)]^{\circ k}(F)\in {\rm Ob}\,\F^{\omega -cons
    (i-1)}$ par l'hypothèse  $H_{div}(A,i)$. Quitte à augmenter la
  valeur de $k$, on peut supposer $\D[L(i)]^{\circ k}(F)=0$. Le
  lemme~\ref{lm-cld} procure alors la conclusion souhaitée : $F\in {\rm Ob}\,\F^{\omega - cons(i-1)}$.

On suppose maintenant $i\leq j\leq n$ et l'assertion vérifiée pour les
valeurs de $j$ strictement inférieures. Comme l'hypothèse
$H^\omega(n)$ est satisfaite et $j$ inférieur à $n$, il existe un foncteur fini
$X$ de $\F_{\Gr,j}$, dont nous noterons $d$ le degré, et une flèche
$(F:A)\xrightarrow{f}\omega_j(X)$ de noyau et conoyau induits de
hauteur strictement inférieure $j$. Si $X=0$, $(F:A)$ est en fait induit de
hauteur au plus $j-1$, et l'hypothèse de récurrence sur $j$ donne la
conclusion. On suppose donc $d\geq 0$, et l'on peut également supposer
que notre assertion est vérifiée pour les valeurs strictement
inférieures de $d$. Soient $u : F\to \D[A](F)\otimes A$ l'unité de l'adjonction, $N$ son
noyau et $Q$ son image. Sous-foncteur de $\D[A](F)\otimes A$, $Q$ est
induit de hauteur au plus $j$ (grâce à l'hypothèse
$H^\omega(n)$). On en déduit que $Q$ est de présentation finie, donc
$N$ de type fini. De plus,
l'application du foncteur $\D[A]$ à la suite exacte
\begin{equation}\label{eq-ser}0\to N\to F\to Q\to 0
\end{equation}
procure une suite exacte
$$\D^1[A](Q)\to\D[A](N)\to\D[A](F)\to\D[A](Q)\to 0$$
dont la flèche $\D[A](F)\to\D[A](Q)$ est un isomorphisme, par
construction de $Q$ (laquelle assure que $\D[A](u)$ est injectif), de sorte que $\D[A](N)$ est un quotient de
$\D^1[A](Q)$. On a par conséquent
$||\D[A](N)||<||Q||=||\D[A](Q)||=||\D[A](F)||=(j,d)$ par la
proposition~\ref{pr-tdgrf} --- en effet, $A(\FF^j)\neq 0$ car $j\geq
i$ et $H_{div}(A,i)$ est satisfaite (ce qui implique $A(\FF^i)\neq 0$ --- cf. remarque~\ref{rqft}).

L'hypothèse de récurrence (sur $j$ si la première composante de
$||\D[A](N)||$ est strictement inférieure à~$j$, sur $d$ sinon)
s'applique donc à $N$ (nous avons déjà remarqué que ce foncteur est de
type fini). Elle montre que $N\in {\rm Ob}\,\F^{\omega-cons(j)}$. La
suite exacte $(\ref{eq-ser})$ établit ensuite $F\in {\rm
  Ob}\,\F^{\omega-cons(j)}$, d'où la proposition.
\end{proof}

\section{Structure du foncteur $P^{\otimes 3}\otimes F$ pour $F$ fini}

La section précédente fournit les outils nécessaires pour
manier des propriétés d'annulation de foncteurs de division de type
$H_{div}$. Dans cette section, nous donnons d'abord des énoncés
généraux pour vérifier de nouvelles propriétés de ce type
à partir d'hypothèses relatives aux simples de $\F$ ou aux représentations des groupes
linéaires (§\,\ref{par31}). Les deux paragraphes suivants les
appliquent pour démontrer d'abord $H_{div}(2)$ (§\,\ref{par32}), ce
qui permet d'affiner des résultats déjà connus, puis d'établir
$H_{div}(3)$  (§\,\ref{par33}), avec pour conséquence le théorème de
structure des foncteurs $P^{\otimes 3}\otimes F$ (pour $F\in {\rm Ob}\,\F^f$).

\subsection{Cadre formel}\label{par31}

\begin{hyp} Dans tout ce paragraphe, on se donne un entier naturel $n$ tel que
l'hypothèse $H^\omega(n)$ est vérifiée.
\end{hyp}

\begin{nota}\label{not-test} On désigne par $\T_n$ la classe des foncteurs $T$ de $\F$ à valeurs de
dimension finie ayant la propriété suivante : {\em si $F\in {\rm
    Ob}\,\F^{tf}$ est tel que $(F:T)\in {\rm Ob}\,\F^{\omega-cons(n)}$,
  alors $F\in {\rm Ob}\,\F^{\omega-cons(n)}$.}
\end{nota}

 La lettre $T$ est
utilisée comme abréviation du terme {\em test} : la division par $T$ teste
l'appartenance à la classe $\F^{\omega-cons(n)}$. Le but de ce
paragraphe consiste à montrer que si l'on sait montrer que
suffisamment de foncteurs appartiennent à la classe $\T_n$, alors on
peut obtenir des renseignements sur les foncteurs de division par
d'autres foncteurs --- précisément, démontrer la condition
$H_{div}(A,n+1)$ pour certains foncteurs~$A$.

Nous donnons d'abord un critère valable pour un foncteur quelconque.

\begin{pr}\label{pr-hd1} Soit $A$ un foncteur de $\F$ à valeurs de
dimension finie vérifiant la propriété suivante : pour tout entier
$m>n$, il existe $T\in\T_n$ tel que pour tout $(\lambda,\mu)\in\p_m$
on ne puisse avoir simultanémement $(S_\lambda : A)=0$, $(S_\mu :
A)=0$ et $S_\lambda$ facteur de composition de $T\otimes S_\mu$.

Alors l'hypothèse $H_{div}(A,n+1)$ est satisfaite.
\end{pr}

\begin{proof} Soit $F\in {\rm Ob}\,\F^{tf}$ tel que $(F:A)=0$, il
s'agit de montrer que $F$ est induit de hauteur au plus $n$. Posons
$m=\max\{\lambda_1\,|\,\lambda\in\p,\;{\rm hom}_\F (F,S_\lambda)\neq
0\}$. Alors $F$ est quotient d'une somme directe finie $P_F$ de projectifs
indécomposables $P_\lambda$ avec $(S_\lambda : A)=0$ (puisque $\D[A]$
est exact à droite) et $\lambda_1\leq m$
(cf. remarque~\ref{rqut-int}). Comme un tel foncteur $P_\lambda$ est
facteur direct de $P_{\FF^m}$, qui est objet de $\F^{\omega - cons(m)}$,
l'hypothèse $H^\omega(n)$ permet d'en déduire $F\in {\rm
Ob}\,\F^{\omega-cons(n)}$ si $m\leq n$.

On suppose donc $m>n$ et l'assertion vérifiée pour les valeurs de $m$
strictement inférieures. Soit $T\in\T_n$ comme dans l'énoncé ;
on examine le foncteur $(F : T)$. Posons $\p(F,T)=\{\mu\in\p\,|\,{\rm hom}_\F ((F:T),S_\mu)\neq
0\}$  et $m'=\max\{\mu_1\,|\,\mu\in\p(F,T)\}$.

D'une part, comme $(P_F : T)\twoheadrightarrow (F : T)$, la relation $\mu\in\p(F,T)$ entraîne qu'il existe
$\lambda\in\p$, avec $(S_\lambda : A)=0$ et $\lambda_1\leq m$, tel que ${\rm hom}_\F((P_\lambda :
T),S_\mu)\neq 0$. Par adjonction, cela équivaut à
${\rm hom}_\F(P_\lambda,T\otimes S_\mu)\neq 0$, i.e. à $S_\lambda$ facteur de
composition de $T\otimes S_\mu$. Cette condition implique elle-même
$\mu_1\leq\lambda_1\leq m$ (évaluer sur
$\FF^{\lambda_1}$).

D'autre part, comme $((F : T):A)\simeq ((F:A):T)=0$, si
$\mu\in\p(F,T)$, alors $(S_\mu : A)=0$.

En combinant les deux points précédents, on obtient que si
$\mu\in\p(F,T)$ est tel que $\mu_1\geq m$, alors $\mu\in\p_m$, $(S_\mu
: A)=0$, et il existe $\lambda\in\p_m$ tel que $(S_\lambda : A)=0$ et
que $S_\lambda$ soit facteur de composition de $T\otimes S_\mu$. Ces
conditions étant incompatibles par hypothèse, il vient
$m'<m$. L'hypothèse de récurrence montre alors que $(F :T)$, qui est
annihilé par $\D[A]$, est induit de hauteur au plus $n$. Comme
$T\in\T_n$, on en déduit que $F$ est également induit de hauteur
au plus $n$, ce qui achève la démonstration.
\end{proof}

L'hypothèse de la proposition précédente est très contraignante. Nous
abordons maintenant un critère plus souple pour la satisfaction d'une hypothèse du
type $H_{div}$ ; en contre-partie, il ne s'applique qu'à des objets
injectifs.

On rappelle que la notation $R_\lambda$ pour les $GL_n$-modules
simples est introduite dans l'appendice.

\begin{pr}\label{pr-dhi} Soient $J$ un objet injectif de co-type
fini  de $\F$ et $\p^{div}(J)=\{\lambda\in\p\,|\,(S_\lambda
: J)=0\}$. On pose également, pour $m\in\mathbb{N}$,
$\p^{div}_m(J)=\p^{div}(J)\cap\p_m$. On suppose que l'hypothèse
suivante est vérifiée : pour tout entier $m>n$, il existe un ordre
total sur $\p^{div}_m(J)$ tel que, pour tout $\lambda\in\p^{div}_m(J)$,
il existe $T_\lambda\in\T_n$ tel que si $R_\mu$ (où $\mu\in\p^{div}_m(J)$) est
facteur de composition du $GL_m$-module 
$R_\lambda\otimes T_\lambda(\FF^m)^*$, alors
$\mu<\lambda$.

Alors  l'hypothèse $H_{div}(J,n+1)$ est satisfaite.
\end{pr}

\begin{proof} Comme le foncteur $\D[J]$ est un facteur direct d'une
somme directe finie de foncteurs décalages (car $J$ est facteur direct
d'une somme directe finie de foncteurs $I_V$), tout foncteur $F$ de $\F$ possède un plus grand quotient
annihilé par $\D[J]$, que nous noterons $\Q_J(F)$ --- précisément, le
foncteur $\D[J]$ possède un adjoint à droite du type $-\otimes G$
(cf. \cite{GP4}, appendice), où
$G$ est un foncteur projectif de type fini, $\Q_J(F)$ est le
conoyau de la coünité $\Q_J(F)\otimes G\to F$. On obtient ainsi un
endofoncteur exact à droite $\Q_J$ de $\F$.

Il s'agit de montrer que pour tout $F\in {\rm Ob}\,\F^{tf}$, on a
$\Q_J(F)\in {\rm Ob}\,\F^{\omega - cons(n)}$. Par exactitude à droite de
$\Q_J$, il suffit de le faire lorsque $F$ est {\em projectif} de type
fini, on encore lorsque $F$ est un foncteur de Powell
$Q_\lambda$ (cf. définition~\ref{df-pow}), puisque les projectifs de type
fini de $\F$ possèdent une filtration finie dont les sous-quotients
sont des foncteurs de Powell.

On procède par récurrence sur
$m=\lambda_1$. Pour $m\leq n$, le foncteur $Q_\lambda$ est lui-même
objet de $\F^{\omega - cons(n)}$, donc l'hypothèse $H^\omega(n)$
entraîne le résultat. On suppose donc $m>n$ et établi que $\Q_J(Q_\mu)\in {\rm
Ob}\,\F^{\omega - cons(n)}$ pour $\mu\in\p_i$ avec $i<m$.

On remarque d'abord que si $\lambda\notin\p^{div}_m(J)$, alors le cosocle $S_\lambda$ de
$Q_\lambda$ (cf. proposition \ref{prpw}) n'est pas annulé par $\D[J]$, donc aucun quotient de $Q_\lambda$
n'est annulé par $\D[J]$, soit $\Q_J(Q_\lambda)=0$. On suppose donc
$\lambda\in\p^{div}_m(J)$.  Une  récurrence sur l'ordre total de
l'énoncé permet de supposer que $\Q_J(Q_\mu)\in {\rm
Ob}\,\F^{\omega - cons(n)}$ pour $\mu\in\p_m^{div}(J)$ tel que $\mu<\lambda$.

Notons $X=(\Prl_m\rho_m(R_\lambda):\iota(T_\lambda))$ : d'après les résultats du
paragraphe~\ref{par-dgr}, $X$ est un foncteur fini
pseudo-constant de $\F_\Gr$, on a $\R_i(X)=0$ pour $i>m$,
$\R_m(X)\simeq\rho_m(R_\lambda\otimes T_\lambda(\FF^m)^*)$ et $(Q_\lambda :
T_\lambda)\simeq\omega(X)$. Par conséquent, il existe une suite exacte
$$0\to F\to(Q_\lambda :
T_\lambda)\to\omega_m\rho_m(R_\lambda\otimes T_\lambda(\FF^m)^*)\to
0$$ où
$F=\omega\Prl_{\leq m-1}\R_{\leq m-1}(X)$, foncteur qui possède une filtration finie
dont les sous-quotients sont des foncteurs de Powell du type $Q_\mu$
avec $\mu_1<m$.

\'Etablissons maintenant que $\Q_J(Q_\lambda : T_\lambda)\in {\rm
Ob}\,\F^{\omega - cons(n)}$.  L'hypothèse de récurrence sur $m$ montre déjà que $\Q_J(F)\in {\rm
Ob}\,\F^{\omega - cons(n)}$. 

Pour prouver que $\Q_J(\omega_m\rho_m(R_\lambda\otimes T_\lambda(\FF^m)^*))\in {\rm
Ob}\,\F^{\omega - cons(n)}$, on note que le foncteur\linebreak[4] $\omega_m\rho_m(R_\lambda\otimes
T_\lambda(\FF^m)^*)$ s'obtient par extensions successives des foncteurs
$Q_\mu$, où $\mu\in\p_m$ est tel que $R_\mu$ est facteur de
composition de $R_\lambda\otimes T_\lambda(\FF^m)^*$. Il suffit donc de
montrer que $\Q_J(Q_\mu)\in {\rm
Ob}\,\F^{\omega - cons(n)}$ pour ces partitions $\mu$. On distingue
pour cela deux cas :
\begin{enumerate}\item $\mu\in\p^{div}_m(J)$. On a donc $\mu<\lambda$,
par choix de $T_\lambda$. L'hypothèse de
récurrence sur $\p^{div}_m(J)$ donne alors $\Q_J(Q_\mu)\in {\rm
Ob}\,\F^{\omega - cons(n)}$.
\item Si $\mu\notin\p^{div}_m(J)$, alors $\Q_J(Q_\mu)=0$ (cf. supra).
\end{enumerate}

On a donc $\Q_J(\omega_m\rho_m(R_\lambda\otimes T_\lambda(\FF^m)^*))\in {\rm
Ob}\,\F^{\omega - cons(n)}$, puis $\Q_J(Q_\lambda : T_\lambda)\in {\rm
Ob}\,\F^{\omega - cons(n)}$. Le foncteur $(\Q_J(Q_\lambda):T_\lambda)$ est un
quotient de $(Q_\lambda : T_\lambda)$ annulé par $\D[J]$, car $\D[T_\lambda]$ et
$\D[J]$ commutent, c'est donc un quotient de $\Q_J(Q_\lambda : T_\lambda)$,
d'où $(\Q_J(Q_\lambda):T_\lambda)\in {\rm
Ob}\,\F^{\omega - cons(n)}$ puisque $H^\omega(n)$ est vérifiée. \'Etant donné que $T_\lambda\in\T_n$, cela entraîne $\Q_J(Q_\lambda)\in {\rm
Ob}\,\F^{\omega - cons(n)}$ comme souhaité.
\end{proof}

\subsection{La structure de $P^{\otimes 2}\otimes F$ pour $F$ fini revisitée}\label{par32}

Le théorème de simplicité généralisé (théorème~\ref{tsg}) et la
proposition~\ref{rap-f} impliquent formellement que l'hypothèse
$H^\omega(1)$ est vérifiée. \`A partir de cette seule observation,
des résultats de la section~\ref{s-div}, utilisée via la proposition suivante, et de la
proposition~\ref{pr-hd1}, nous allons établir l'importante proposition~\ref{pr-hf2}.

\begin{pr}\label{pr-fivl1} L'hypothèse $H_{div}(\Lambda^1,1)$ est satisfaite.
\end{pr}

\begin{proof} Si $F$ est un foncteur fini de degré $n$ de $\F$, alors le foncteur
  $(F:\Lambda^1)$ est de degré $n-1$, donc non nul, si $F$ est non
  constant (cf. \cite{art1}, corollaire 3.10). Comme $\D[\Lambda^1]$ est exact à droite, on en déduit
  qu'un foncteur $F$ tel que $(F:\Lambda^1)=0$ n'a des quotients finis
  que constants. Comme $F$ est limite de ses quotients finis
  (cf. proposition~\ref{rap-f}), cela montre que $F$ est constant,
  donc fini, $F$ étant de type fini, d'où la conclusion.
\end{proof}

\begin{pr}\label{pr-hf2} L'hypothèse $H_{div}(S_{<2>},2)$ est satisfaite.
\end{pr}

\begin{proof} Nous avons déjà rappelé
  que $H^\omega(1)$ est vérifiée. Comme l'hypothèse $H_{div}(1)$ est
  vérifiée ($\D[L(1)]\simeq\Delta$ n'annihile que les foncteurs constants), les propositions~\ref{pr-fivl1}
  et~\ref{pr-omd} impliquent qu'un foncteur $F\in {\rm Ob}\,\F^{tf}$
  tel que $(F : \Lambda^1)\in {\rm Ob}\,\F^{\omega - cons(1)}$ est
  lui-même objet de $\F^{\omega - cons(1)}$. Autrement dit, avec la
  terminologie du paragraphe~\ref{par31}, on a $\Lambda^1\in\T_1$ ; c'est le seul
  foncteur~\guillemotleft~test~\guillemotright~dont nous aurons besoin
  ici. Vérifions à présent que l'hypothèse de la
  proposition~\ref{pr-hd1} est satisfaite pour $A=S_{<2>}$  : la proposition~\ref{pr-divst}
  montre que les partitions régulières $\lambda$ telles que
  $(S_\lambda : S_{<2>})=0$ vérifient $l(\lambda)\leq 1$. La
  conclusion résulte alors de ce que $\Lambda^m$ n'est pas facteur de
  composition de $\Lambda^1\otimes\Lambda^m$ si
  $m>1$.
\end{proof}

On retrouve en particulier la satisfaction de $H^\omega(2)$, déduite
dans \cite{artsmf} (§\,16.4) des résultats de \cite{GP5} sur le foncteur
$\tilde{\nabla}_2$ :

\begin{cor}\label{crf1} Les hypothèses $H_{div}(2)$ et $H^\omega(2)$
sont vérifiées.
\end{cor}

\begin{proof} Comme $S_{<2>}$ est un sous-foncteur de $L(2)$, $\D[S_{<2>}]$
est un quotient de $\D[L(2)]$, de sorte que $H_{div}(S_{<2>},2)$
implique $H_{div}(2)$. La proposition~\ref{pr-cadiv} montre alors
qu'il en est de même pour $H^\omega(2)$.
\end{proof}

On retrouve ainsi le caractère noethérien de type~$2$ de $P^{\otimes
2}\otimes F$ pour $F\in {\rm Ob}\,\F^f$.

\begin{rem} Les renseignements obtenus par notre démarche sont un
petit peu plus précis que ceux donnés par l'approche de Powell dans
\cite{GP5}. En effet, par le corollaire~\ref{cr-andiv}, on voit que le
corollaire~\ref{crf1} décrit tous les foncteurs de type fini $F$
n'ayant pour facteurs de composition que des puissances extérieures
(modulo les foncteurs finis, ce sont les quotients des sommes directes
d'un nombre fini de copies de $P_\FF$). En revanche, l'article
\cite{GP5} utilise fortement que les foncteurs de type fini $F$ tels
que $\tilde{\nabla}_2(F)=0$ non seulement n'ont que des puissances
extérieures comme facteurs de composition, mais n'ont pas non plus de
sous-quotients isomorphes à certaines extensions entre puissances
extérieures. De plus, on constate que \cite{GP5} nécessite la
connaissance de ${\rm Ext}^1_\F(\Lambda^i,\Lambda^j)$ (groupes
déterminés par Franjou dans \cite{Franjou}), alors que notre
démonstration n'utilise, comme point non formel, que le fait très
élémentaire que $\Lambda^i$ n'est facteur de composition de
$\Lambda^1\otimes\Lambda^i$ que pour $i=1$.
\end{rem}

Nous terminons ce paragraphe par une autre conséquence directe de la
proposition~\ref{pr-hf2}, fondée sur la remarque~\ref{rqft}.\,\ref{itf}.

\begin{pr}\label{pr-ldx} L'hypothèse $H_{div}(\Lambda^2,2)$ est satisfaite.
\end{pr}

\begin{proof} Si $(F : \Lambda^2)=0$, alors $(F :
  \Lambda^2\otimes\Lambda^1)=0$ ; comme $S_{<2>}$ est facteur direct
  de $\Lambda^2\otimes\Lambda^1$, on en déduit $(F : S_{<2>})=0$. La
  conclusion résulte donc de la proposition~\ref{pr-hf2}. 
\end{proof}

\subsection{Les résultats principaux}\label{par33}

\begin{theo}\label{th-pr} L'hypothèse $H_{div}(3)$ est vérifiée.
\end{theo}

\begin{proof} Nous déduirons ce résultat de la
  proposition~\ref{pr-dhi}. Comme $H^\omega(2)$ est vérifiée, la
  satisfaction de $H_{div}(2)$, $H_{div}(S_{<2>},2)$,
  $H_{div}(\Lambda^1,1)$ et $H_{div}(\Lambda^2,2)$ (cf. paragraphe
  précédent) implique, par la proposition~\ref{pr-omd}, que $S_{<2>}$,
  $\Lambda^1$ et $\Lambda^2$ sont dans $\T_2$
  (cf. notation~\ref{not-test}). 

La proposition~\ref{pr-divst} nous apprend que
$\p^{div}(L(3))=\{\lambda\in\p\,|\,l(\lambda)\leq 2\}$, selon la
notation de la proposition~\ref{pr-dhi}. Pour $m>2$, on munit
l'ensemble $\p^{div}_m(L(3))$ de l'ordre total donné comme suit :
$$(m)<(m,1)<(m,2)<\dots<(m,m-1).$$

Comme $R_{(m)}$ est l'unité du produit tensoriel des $GL_m$-modules,
on a $$R_{(m)}\otimes S_{<2>}(\FF^m)^*\simeq R_{(m,2,1)}^*\simeq
R_{(m,m-1,m-2)}.$$ Ce $GL_m$-module n'a pas de facteur de composition
$R_\mu$ avec $\mu\in\p^{div}_m(L(3))$.

Pour $0<i<m-1$, on a $$R_{(m,i)}\otimes\Lambda^1(\FF^m)^*\simeq
R_{(m,i)}\otimes R_{(m,1)}^*\simeq R_{(m,i)}\otimes R_{(m,m-1)}\simeq
(\Lambda^i\otimes\Lambda^{m-1})(\FF^m).$$
Comme $i<m-1$, les facteurs de composition du foncteur
$\Lambda^i\otimes\Lambda^{m-1}$ sont les $S_{(m-1+t,i-t)}$ pour $0\leq
t\leq i$, donc les facteurs de composition de
$(\Lambda^i\otimes\Lambda^{m-1})(\FF^m)$ sont $R_{(m,m-1,i)}$ et
$R_{(m,i-1)}$. Ainsi, $\mu\in\p^{div}_m(L(3))$ et $R_\mu$ facteur de
composition de $R_{(m,i)}\otimes\Lambda^1(\FF^m)^*$ implique
$\mu=(m,i-1)<(m,i)$ (ou $\mu=(m)<(m,1)$ si $i=1$).

De façon analogue, $R_{(m,m-1)}\otimes\Lambda^2(\FF^m)^*\simeq
(\Lambda^{m-1}\otimes\Lambda^{m-2})(\FF^m)$ a pour facteurs de
composition $R_{(m,m-1,m-2)}$ et $R_{(m,m-3)}$ ($R_{(3)}$ si
$m=3$). Donc si $\mu\in\p^{div}_m(L(3))$ est tel que $R_\mu$ est facteur de
composition de $R_{(m,i)}\otimes\Lambda^1(\FF^m)^*$, on a
$\mu=(m,m-3)<(m,m-1)$ (ou $\mu=(3)<(3,2)$ si $m=3$).

La proposition~\ref{pr-divst} procure alors la satisfaction de $H_{div}(L(3),3)=H_{div}(3)$.
\end{proof}

Ce théorème donne un renseignement bornant la~\guillemotleft~taille~\guillemotright~des foncteurs de type fini
dont les facteurs de composition sont limités. Précisément : 

\begin{cor}\label{crfc} Un foncteur de type fini de $\F$ dont tous les facteurs de
composition sont associés à des partitions de longueur au plus $2$ est
induit de hauteur au plus~$2$.
\end{cor}

\begin{proof} Par le
  corollaire~\ref{cr-andiv}, les foncteurs annihilés par $\D[L(3)]$
  sont exactement ceux dont les facteurs de composition sont du type
  $S_\lambda$ avec $l(\lambda)\leq 2$. Le théorème~\ref{th-pr} donne
  donc la conclusion.
\end{proof}

\begin{rem} Il semble très difficile de donner des informations
  beaucoup plus explicites sur les foncteurs de type fini dont les facteurs de
composition sont associés à des partitions de longueur au plus
$2$. Cette observation suggère qu'on ne peut pas trouver d'approche
directe pour démontrer la corollaire~\ref{crfc}.
\end{rem}

L'avancée principale sur la structure globale de la catégorie~$\F$
apportée par le théorème~\ref{th-pr} est le résultat suivant.

\begin{cor} L'hypothèse $H^\omega(3)$ est vérifiée : la sous-catégorie
$\F^{\omega - cons(3)}$ des foncteurs induits de hauteur au plus
$3$ de $\F$ est épaisse.
\end{cor}

\begin{proof} Ce résultat se déduit du théorème~\ref{th-pr} et de la proposition~\ref{pr-cadiv}.
\end{proof}

Compte-tenu de la proposition~\ref{cortsg}, on en déduit le
résultat suivant, qui semble constituer la meilleure forme partielle
démontrée de la
conjecture artinienne.

\begin{cor} Pour tout foncteur fini $F$ de $\F$, le foncteur
$P^{\otimes 3}\otimes F$ est noethérien de type~$3$.
\end{cor}

\appendix

\section{La catégorie $\F$ et les représentations des groupes linéaires}

Les propriétés élémentaires de la catégorie $\F$ rappelées dans cet appendice, qui sert surtout à fixer les notations, sont exposées et démontrées en détails dans les articles fondamentaux \cite{K1} et \cite{K2} (où l'on trouvera également d'autres références). 

\paragraph*{Propriétés de finitude dans une catégorie de Grothendieck
  $\A$ (cf. \cite{Pop}).} Les objets {\em simples} de $\A$ sont les objets non nuls n'ayant pas
de sous-objet strict non nul. Les objets de {\em longueur finie} (i.e. possédant une filtration
finie de sous-quotients simples) sont simplement appelés objets {\em
  finis} dans cet article ; la sous-catégorie pleine des objets finis de
$\A$ est notée $\A^f$. Les objets localement finis sont les colimites d'objets
finis de $\A$ ; ils forment une sous-catégorie pleine $\A^{lf}$. 

La {\em filtration de Krull} de $\A$ est la suite croissante  $(\K_n(\A))$ de
sous-catégories localisantes (c'est-à-dire épaisses et stables par
colimites) de $\A$ définie par récurrence par
$\K_n(\A)=\{0\}$ si $n<0$ et le fait que, pour $n\in\mathbb{N}$,
$\K_n(\A)/\K_{n-1}(\A)$ est la sous-catégorie localisante de
$\A/\K_{n-1}(\A)$ engendrée par les objets finis de cette catégorie. Un objet de
$\A$ est dit {\em noethérien de type $n$} s'il est noethérien et
appartient à $\K_n(\A)$ (la définition usuelle d'objet noethérien de type $n$ et
son équivalence avec celle-là sont rappelées à la fin de l'appendice B
de \cite{artsmf}). La sous-catégorie pleine des objets
noethériens de type~$n$ de $\A$ sera notée $\mathbf{NT}_n(\A)$ ; elle
est épaisse.

Les objets {\em de type fini} de $\A$ sont les
objets $A$ tels que toute  suite croissante de sous-objets de $A$ de
colimite $A$ stationne. Ainsi, les objets noethériens sont ceux dont
tous les sous-objets sont de type fini. Dans la catégorie $\F$, les objets de
type fini sont exactement les quotients des sommes directes finies de
foncteurs projectifs standard. La sous-catégorie pleine des objets de
type fini de $\A$ sera notée $\A^{tf}$. 

Un objet $A$ de $\F$ est de {\em présentation finie} s'il existe une suite
exacte $P'\to P\to A\to 0$ avec $P$ et $P'$ projectifs de type
fini (cette définition est adaptée car $\F$ a
  assez de projectifs de type fini).

Dans $\F$, nous dirons qu'un objet est de {\em co-type fini} si son
dual est de type fini.

\paragraph*{Les simples de $\F$.} Les objets simples de $\F$ peuvent se classifier à l'aide des
représentations irréductibles (sur $\FF$) des groupes
symétriques. On utilise pour cela les {\em partitions}, i.e. les suites finies
décroissantes d'entiers strictement positifs $\lambda=(\lambda_1,\dots,\lambda_r)$
($r$ s'appelle alors la {\em longueur} de $\lambda$ et se note
$l(\lambda)$ ; la somme des $\lambda_i$ s'appelle le {\em degré} de
$\lambda$ et se note $|\lambda|$). Par convention, on posera
$\lambda_i=0$ pour $i>l(\lambda)$.

Si $\lambda$ est une partition, notons
$\Lambda^\lambda$ le produit tensoriel de puissances extérieures
$\Lambda^{\lambda_1}\otimes\dots\otimes\Lambda^{\lambda_r}$. Le
foncteur $\Lambda^\lambda$ possède un sous-foncteur remarquable appelé
{\em foncteur de Weyl} associé à $\lambda$ et noté~$W_\lambda$. 

\begin{rem}\label{rq-gst} Les éléments {\em semi-standard} de $\Lambda^\lambda(V)$, c'est-à-dire les éléments du type
$$(u_1\wedge\dots\wedge u_{\lambda_1})\otimes(u_1\wedge\dots\wedge u_{\lambda_2})\otimes\dots\otimes
  (u_1\wedge\dots\wedge u_{\lambda_{l(\lambda)}}),$$
où $u_1,\dots,u_{\lambda_1}$ sont des éléments de $V$, appartiennent à $W_\lambda(V)$.

Lorsque la partition $\lambda$ est régulière (notion rappelée ci-dessous), le foncteur $W_\lambda$ est engendré par les éléments semi-standard. Cela résulte de la
  variante fonctorielle des résultats classiques de James \cite{James} qui est présentée par exemple dans l'article
  \cite{PS1}\,\footnote{On prendra garde que nos conventions en matière
  de partitions sont duales de celles de \cite{PS1}.}.
\end{rem}

Si $\lambda$ est une partition
{\em $2$-régulière} (nous dirons simplement {\em régulière} dans cet article, le corps de base étant fixé à $\FF$), i.e. une suite finie {\em strictement} décroissante
d'entiers strictement positifs, le
cosocle (i.e. plus grand quotient semi-simple) de $W_\lambda$ est un foncteur simple noté $S_\lambda$. Lorsque $\lambda$ parcourt l'ensemble des partitions régulières,
$S_{\lambda}$ décrit un système complet de représentants des objets
simples de $\F$ modulo isomorphisme (cf. \cite{PS1}).

Les objets simples de $\F$ sont {\em auto-duaux} : $DS_\lambda\simeq
S_\lambda$ (l'auto-dualité est en fait une condition un petit peu plus
forte --- cf. \cite{PS1}).

\begin{nota}\label{not-part}\begin{itemize}\item On désigne par $\p$
    l'ensemble des partitions régulières.
\item Pour $i\in\mathbb{N}$, on pose $\p_i=\{\lambda\in\p\,|\,\lambda_1=i\}$.
\item Pour $\lambda\in\p$, on note $P_\lambda$ la couverture
  projective de $S_\lambda$.
\item Pour $n\in\mathbb{N}$, on note $<n>$ la partition $(n,n-1,\dots,1)$ de $n(n+1)/2$.
\end{itemize}
\end{nota}

\begin{ex}\begin{enumerate}\item Les foncteurs puissances extérieures
    sont simples : $S_{(i)}=\Lambda^i$.
\item Pour $i>j>0$, les facteurs de composition (i.e. sous-quotients
  simples) de
  $\Lambda^i\otimes\Lambda^j$ sont exactement les $S_{(i+t,j-t)}$ pour
  $0\leq t\leq j$ (en convenant que $S_{(k,0)}=\Lambda^k$). Les
  facteurs de composition de $\Lambda^i\otimes\Lambda^i$ sont
  $\Lambda^i$ et les $S_{(i+t,i-t)}$ pour $0<t\leq i$.
\item Le foncteur $S_{<2>}$ est caractérisé par le scindement
  $\Lambda^2\otimes\Lambda^1\simeq\Lambda^3\oplus S_{<2>}$.
\end{enumerate}
\end{ex}

\paragraph*{Les $GL_n$ - modules simples.} On note $GL_n=GL_n(\FF)$
le groupe linéaire sur $\FF$ et $GL_n -
\mathbf{mod}$ la catégorie des $\FF[GL_n]$-modules à gauche, appelés
simplement $GL_n$-modules dans cet article.

Pour $\lambda\in\p_n$, le $GL_n$-module
$S_\lambda(\FF^n)$, que nous noterons $R_\lambda$, est simple, et
les $R_\lambda$ constituent un système complet de représentants des
$GL_n$-modules simples lorsque $\lambda$ décrit $\p_n$ (cf. \cite{K2} ; on prendra garde que les
conventions d'indexation de cet article de Kuhn, de l'article
\cite{PS1} de Piriou et Schwartz et du présent travail sont deux à deux
distinctes).

Pour $\lambda\in\p$ et $n\in\mathbb{N}$ tel que $\lambda\notin\p_n$
(cas traité précédemment), on a $S_\lambda(\FF^n)=0$ si
$n<\lambda_1$ et $S_\lambda(\FF^n)\simeq R_{(n,\lambda)}$ si
$\lambda_1>n$.

Nous noterons $(-)^*$ le foncteur de dualité contragrédiente de $GL_n -
\mathbf{mod}$. Son effet sur les simples est donné par
$R_{\lambda}^*\simeq
R_{(n,n-\lambda_r,\dots,n-\lambda_2)}$ pour $\lambda\in\p_n$, où
$r=l(\lambda)$. Cette observation, qui se déduit aussitôt de
\cite{K2} (§\,6) et de ce que les représentations $\Lambda^i(\FF^n)$ et
$\Lambda^{n-i}(\FF^n)$ de $GL_n$ sont duales, est utilisée dans la
démonstration du théorème~\ref{th-pr}.

\begin{rem}\label{rqstein} La partition régulière $<n>$ joue un rôle important dans cet article ;
il convient de noter que $R_{<n>}$ est la {\em représentation de
  Steinberg} de $GL_n(\FF)$, assertion qui a déjà
été exploitée dans la section~13.1 de \cite{artsmf}. Une des propriétés utiles est que la projection
canonique $W_{<n>}\twoheadrightarrow S_{<n>}$ est un isomorphisme
--- cf. \cite{James}, corollaire $24.9$ ; nous renvoyons aussi à
l'article \cite{Mitch} de Mitchell, où l'on
trouvera une démonstration de l'équivalence entre les définitions
classiques de la représentation de Steinberg et la
définition~\guillemotleft~fonctorielle~\guillemotright~$R_{<n>}$. 

L'enveloppe injective de $S_{<n>}$ (dans $\F$) se note
traditionnellement $L(n)$. L'exemple 3.0.3.2 de l'article \cite{GP2} montre que
l'inclusion $S_{<n>}\hookrightarrow L(n)$ est un isomorphisme
lorsqu'on l'évalue sur un espace de dimension au plus $n$ ; en particulier, on a $L(n)(\FF^{n-1})=0$.
\end{rem}

\paragraph*{Foncteur différence et foncteurs finis.}\label{fdif} Pour $E\in {\rm
  Ob}\,\E^f$, on introduit l'endofoncteur $\Delta_E$ de $\F$, appelé {\em décalage par $E$}, donné par $\Delta_E(F)(V)=F(V\oplus E)$.

Le {\em
  foncteur différence} $\Delta$ est caractérisé par le scindement $\Delta_\FF\simeq id\oplus\Delta$. Il commute à toutes les limites et colimites.

Les foncteurs $F$ de $\F$
tels qu'il existe $i\in\mathbb{N}$ tel que $\Delta^{i}(F)=0$ sont
appelés {\em foncteurs polynomiaux}. Le degré d'un foncteur polynomial
$F$ est le plus grand entier positif $\deg F$ tel que $\Delta^{\deg F}
F\neq 0$ si $F\neq 0$ ; par convention, $\deg 0=-\infty$

\begin{ex} Si $\lambda$
est une partition, les foncteurs $\Lambda^\lambda$, et $W_\lambda$ de
$\F$, de même que $S_\lambda$ si $\lambda$ est régulière, sont polynomiaux de degré
$|\lambda|$. 
\end{ex}

Un résultat de base sur la catégorie $\F$ est le suivant : 
\begin{pr}[Cf. \cite{K1}]\label{rap-f}\begin{itemize}\item Un foncteur de $\F$ est
    fini si et seulement s'il est polynomial et à valeurs de dimension finie.
\item Tout foncteur de type fini de $F$ est limite de ses quotients
  polynomiaux (donc finis).
\end{itemize}
\end{pr}

\begin{rem}\label{rqut-int}\begin{itemize}\item Une conséquence utile de la deuxième assertion est la suivante : pour tout
$F\in {\rm Ob}\,\F^{tf}$, l'épimorphisme canonique
$F\twoheadrightarrow {\rm cosoc}\,F$ de $F$ sur son cosocle est essentiel, donc, si l'on écrit
${\rm cosoc}\,F\simeq\bigoplus_{\lambda\in\p}S_\lambda^{\oplus a_\lambda}$, alors $F$ est
isomorphe à un quotient de $\bigoplus_{\lambda\in\p}P_\lambda^{\oplus
  a_\lambda}$, qui en est une couverture projective.
\item Un corollaire de la première assertion est la stabilité par
  produit tensoriel des
  objets finis dans la catégorie~$\F$.
\end{itemize}
\end{rem}

\nocite{*}
\bibliographystyle{smfalpha}
\bibliography{bibca3}
\end{document}